\documentclass[12pt]{amsart}
\usepackage{amssymb, latexsym}
\usepackage{latexsym}
\usepackage{amsmath}
\usepackage{amsfonts}
\theoremstyle{plain}
\newtheorem*{theorem}{Theorem}
\newtheorem{corollary}{Corollary}

\newtheorem{lemma}{Lemma}
\newtheorem{proposition}{Proposition}
\theoremstyle{definition}

\theoremstyle{remark}

\newtheorem{remark}{Remark}
\numberwithin{equation}{section}
\newcommand{\HN}{\mathcal H_N}

\def\d{\delta}
\def\a{\alpha}
\def\b{\beta}
\def\g{\gamma}

\def\e{\epsilon}
\def\l{\lambda}
\def\m{\mu}

\def\r{\rho}
\def\s{\sigma}
\def\ta{\tau}
\def\th{\theta}

\def\th{\theta}
\def\G{\Gamma}

\def\R{\mathbb R}
\def\C{\mathbb C}
\def\tr{\operatorname{tr}}
\def\sgn{\operatorname{sgn}}
\usepackage{latexsym}
\def\d{\delta}
\def\a{\alpha}
\def\b{\beta}
\def\g{\gamma}

\def\e{\epsilon}
\def\l{\lambda}
\def\m{\mu}

\def\r{\rho}
\def\s{\sigma}
\def\ta{\tau}
\def\th{\theta}

\def\th{\theta}
\def\G{\Gamma}

\def\C{\mathbb C}
\def\N{\mathbb N}
\def\R{\mathbb R}
\def\tr{\operatorname{tr}}
\def\sgn{\operatorname{sgn}}
\def\max{\operatorname{max}}
\def\HN{\mathcal H_N}
\begin{document}
\title[Fixed trace random matrices]
{Limit correlation functions \\
for fixed  trace \\ random matrix ensembles
\date{January 25, 2004}}
\author [Friedrich G\"otze \and Mikhail Gordin ]
 {Friedrich G\"otze$\,^1$ \and Mikhail Gordin$\,^{1,2}$ }
\thanks {$^1$ Research supported
by Sonderforschungsbereich 701 "Spektrale Strukturen und Topologische Methoden in der Mathematik".}
\thanks{$^2$ Partially supported by grants RFBR-05-01-00911, DFG-RFBR-04-01-04000, and NS-2258.2003.1} 
\address
{(F. G\"otze)\\
 Fakult\"at f\"ur Mathematik\\
 Universit\"at Bielefeld\\
Postfach 100131, D-33501, Bielefeld\\
 GERMANY}
\email{goetze@mathematik.uni-bielefeld.de}
\address
{(M. Gordin)\\
 V.A. Steklov Mathematical Institute at St. Petersburg\\
27 Fontanka, St. Petersburg, 191023\\
 RUSSIA}
\email{gordin@pdmi.ras.ru} 
\maketitle
\begin{abstract}
Universal limits for the eigenvalue correlation functions
in the bulk of the spectrum are shown for a class of
nondeterminantal random matrices known as the fixed trace or
the Hilbert-Schmidt ensemble. These universal limits have been
proved before for determinantal Hermitian matrix ensembles and
for some special classes of the Wigner random matrices.
\end{abstract}
\section{Introduction and the statement of the result}
\label{1}
 Let $\HN$ be the set of all $ N\times N $ (complex) Hermitian
matrices, and let $\tr A=\sum_{i=1}^N a_{ii}$ denote the
trace of a square matrix $A=(a_{ij})_{i,j=1}^N.$   $\HN$ is
a real Hilbert space of
dimension $N^2$ with respect to the symmetric bilinear form
 $(A,B) \mapsto \tr AB.$ Let $l_N$ denote the unique Lebesgue measure
 on $\HN$ which satisfies the relation  $l_N(Q)=1$  for every
cube $Q\subset \HN$ with edges of length 1. A Gaussian
probability measure on $\HN$ invariant with respect to all
orthogonal linear transformations of $\HN$ is uniquely
 defined up to a scaling transformation.  Such  measures form a
one-parameter family $(\m_N^{s})_{s>0},$ where the measure
$\m_N^{s}$ is specified by its density
 \begin{equation}
g_N^{s}(A)= \frac{1}{(\sqrt{s2\pi})^{N^2}}
\exp{\biggl(-\frac{1}{2s}  \  \tr A^2\biggr)}
\end{equation}
with respect to $l_N.$ Thus, for a random matrix $X$
distributed according to $\m_N^{s}$ we have
\begin{equation}
E_{\m_N^{s}} \tr X^2= s N^2.
\end{equation}
The set $\HN$ endowed with the measure $\m_N^s$
is called the Gaussian Unitary Ensemble (GUE).
Let $X$ be a random $ N \times N $  Hermitian matrix (that
 is a random variable taking values in $\HN$).  We consider the eigenvalues
 $\l_1, \l_2,  \dots, \l_N$
of the random matrix $X$ as a finite sequence of {\em exchangeable}
random variables.  By definition, this means that their joint distribution
$P_N^X$ does not change under any permutation of these variables.
  Let for each $ n, 1 \le n \le N, $ $P_{n,N}^X$ denote the joint distribution of
  some $n$ of these $N$ variables. Obviously, $P_{n,N}^X$  is a permutation invariant probability measure in $\mathbb R^n $.
In particular, the measure
$P_{1,N}^X$ describes the distribution of a single eigenvalue.
By definition, the $n${\it -point correlation measure}
$\r_{n,N}^X$  of a random matrix $X$ is  a non-normalized
measure defined by the relation
\begin{equation} \r_{n,N}^X=\frac{N!}{(N-n)!}P_{n,N}^X.
\end{equation} For a measurable
set $A \subset \R^n$ the quantity $\r_{n,N}^X(A)$ can be
interpreted as the average number of $n-$tuples of eigenvalues
 in the set $A$. If the measure  $ \r_{n,N}^X$  is
absolutely continuous with respect to the Lebesgue measure on
$\R^n$, its Radon-Nikodym derivative $R_{n,N}^X$ is called the
$n${\it -point correlation function} of the random matrix $X$.
In particular, the measure $\r_{1,N}^X$ has total mass $N$.
For a measurable set $E \subset \R^1$, the quantity $\r_{1,N}^X
(E)$ expresses the expected number of the eigenvalues belonging
to $E$. The corresponding density with respect to the Lebesgue
measure in $\R^1,$ if it exists, is called the
{\it eigenvalue density} or the {\it density of states}
(caution: under the same names the normalized versions of the same
measures are considered in the literature as well).
Let $X_N$ be a random matrix with the distribution  $\m_N^s$.
For $n=1, \dots,N$ we set $P_{n,N}^{\text{GUE},
s}=P_{n,N}^{X_N}$ and $\r_{n,N}^{\text{GUE},s} =
\r_{n,N}^{X_N}$. A classical result for the $\text{GUE}$
says that we have
\begin{equation} P_{1,N}^{\text{GUE}, 1/ N}
\underset{N \to \infty}{ \to} W, \end{equation} where the
measures converge  in the weak sense, and
 $W$ is the standard Wigner measure on
$[-2,2]$ defined by the density
\begin{equation} \label{wigner}
w(x)= (2\pi)^{-1}\sqrt{(4-x^2)_+}, \ x \in  \R. \end{equation}
 In terms of the correlation
measures the same relation reads
\begin{equation}\label{1.1}
\frac{1}{N}\r_{1,N}^{\text{GUE}, 1/ N} \underset{ N \to
\infty}{\to} W.
\end{equation}
For the $n-$point correlation measures we have a similar
relation
\begin{equation}\frac{1}{N^n}\r_{n,N}^{\text{GUE}, 1/ N}
\underset{ N \to \infty}{\to} \underset{n  \  \text{times}} {W
\times \cdots \times W},\end{equation} which means that the
eigenvalues become independent in the limit. However, for
$n \ge 2,$ the study of a finer asymptotics near a point from 
the principal diagonal in 
the cube $(-2,2)^n$  shows \cite{J,M}: for every
$u \in (-2,2)$ and $t_1, \dots, t_n \in \R^1$
\begin{align}\label{sd}
 \lim_{N \to \infty}\frac{1}{(Nw(u))^n} R_{n,N}^{\text{GUE},\tfrac{1}{N}}
\bigl(u+(t_1/Nw(u)), &\dots,
  u+(t_n/Nw(u))\bigr) \\
=&\det\biggl(\frac{\sin \pi(t_i-t_j)}{
\pi(t_i-t_j)}\biggr)_{i.j=1}^n. \notag \end{align}
 This limit relation presents a pattern for many other results,
 in particular,  for that of the present paper.
 The right hand side of this relation represents an example
of the  correlation function of a so-called determinantal
(or fermionic) random point process  \cite{S}. In general the
$n-$point correlation function $R_n $ of such a process is given by the formula
$$R(u_1, \dots, u_n)= \det\bigl|K(u_i,u_j)|_{i,j=1}^n, $$
where  $K$ is the kernel of an integral operator on the line, which is
trace-class having been restricted to any finite interval in $\mathbb R,$ and
subject to some further conditions
(see \cite{S} for a detailed exposition). Moreover, in the
asymptotic Hermitian random matrix theory
$(K_N)_{N\ge0}$ are the reproducing kernels of the subspaces of polynomials of degree $\le N-1$ with respect to some weight on the line. In this case we call the
corresponding matrix ensemble {\em determinantal.} The GUE gives an example of such an ensemble. To a large extent the asymptotic study of determinantal ensembles reduces to that of the respective kernels \cite{D99A,D99C}. In particular, the same 
local limit as in the case of GUE is established in \cite{PS,D99A} for two broad classes of determinantal matrix ensembles. It is expected that the $\sin$-kernel limit (first discovered by F. Dyson) is rather common. This is known as the {\em universality conjecture}.
 Outside the class of determinantal Hermitian random matrices only very few results on the asymptotics of the correlation functions are known (see, for instance, \cite{J}, where a mixture of determinantal measures is considered). \\ 
  In the present paper we investigate the following
non-determinantal
  ensemble of Hermitian random matrices. Let
\begin{equation}S_N^r =\{A \in \HN: \tr A^2=r^2\}\end{equation}
be the sphere in $\HN$ of the radius $r > 0$  centered at the
origin. Set $r=\sqrt{s} N.$ The  sphere  $S_N^{\sqrt{s} N}$
carries a unique probability measure $\nu_N^s $
invariant with respect to
all orthogonal linear transformations in the space $\HN$.
We call this measure the {\em  fixed
Hilbert-Schmidt norm
ensemble} (or just HSE) to reserve the term "the fixed trace ensemble"
for more general ones (see \cite{A}). Let $Y_N$ be a random matrix
distributed according to  $\nu_N^s$. We set for $n=1, \dots,N$
$P_{n,N}^{\text{HSE},s}=P_{n,N}^{Y_N}$ and
$\r_{n,N}^{\text{HSE},s}=\r_{n,N}^{Y_N}$.
It is a well known result (see \cite{M,R}) that
\begin{equation}
\label{HSElim} P_{1,N}^{\text{HSE},1/N} \underset{ N \to
\infty}{\to} W,
\end{equation}
like in the case of GUE.  In this paper we prove that the
correlation functions $R_{n,N}^{\text{HSE},1/N}$ of arbitrary
order $n \  (1\le n \le N)$ near every point $u \in
(-2,2),$ have the same determinantal limit with the
kernel $\sin \pi(t _1- t_2)/\pi (t_1 - t_2)$
 as the GUE correlation functions (for $n=1$ the limit equals $1$).
More precisely, we establish in this paper the following result.
\begin{theorem} Let $R_{n,N}^{\nu,1/N}$ be the $n-$point
correlation function of the eigenvalues for a random
matrix uniformly distributed on the sphere  $S_N^{\sqrt{N}}$.
Then for every $u \in (-2,2)$ and $t_1, t_2, \dots, t_n \in \R^1
$
\begin{equation} \label{conv1}
  \frac{1}{(Nw(u))^n}R_{n,N}^{\nu,1/N}\left(u+\tfrac{t_1}{Nw(u)},
 \dots,u+\tfrac{t_n}{Nw(u)} \right)
 \to \det  \biggl( \frac{\sin \pi(t_i-t_j)}
 {\pi(t_i-t_j)}\biggr)_{\mspace{-4.0 mu} i,j=1}^{\mspace{-4.0 mu}n}  \mspace{-16.0 mu} 
 \end{equation}
as $N \to \infty$. For every $\a \in (0,1)$ and $A > 0$
the relation \eqref{conv1} holds uniformly in all $u \in [-2 + \a,
2 - \a],$ and $t_1 \in [-A,A], \dots, t_n \in [-A,A]$.
\end{theorem}

The class of fixed trace matrix ensembles was studied in several recent publications. The authors are indebted to G. Akemann for drawing their attention, 
after the first version of the present work has appeared as a preprint, 
to the paper \cite{A00C}, he published jointly with G. Vernizzi.  
In this paper the local universality is studied for a class of ensembles containing the HSE. 
The authors provide heuristic arguments for 
the universality near zero based on the complex inversion of the Laplace transform. Our proof is based on the complex Laplace inversion as well (more precisely, 
we arrive at the Fourier transform after a sequence of translations of the origin and rescalings depending on $N$). However, establishing bounds for orthogonal functions and related kernels in the complex domain allows us 
to rigorously conclude convergence near every point of the interval $(-2,2).$

 In the following we will sketch the main steps of the proof. The guiding
principle is that results for
the Hilbert-Schmidt ensembles are deducible from the corresponding
results for the GUE using the 'equivalence of ensembles'
or the 'concentration phenomenon'.
In our setup a primitive form of the concentration is given by the
law of large numbers for the squares of the Hilbert-Schmidt norms
of the GUE random matrices.  Supplemented by some estimates of
the probabilities of large deviations, this is the main
tool in \cite{GGL}, where the universality is shown near zero for the correlation measures (rather than the correlation \linebreak
functions) of the Hilbert-Schmidt ensemble; 
convergence was shown there in the weak topology determined by the continuous compactly supported functions. Such simple arguments seem to be 
 insufficient for proving local results
for correlation functions of eigenvalues (maybe, with exception
for eigenvalues near zero), particularly, for stronger topologies. This agrees with M.L.Mehta's doubts (\cite{M},
Sect. 27.1, p.490) concerning the deducibility of the local results for
correlation functions of the Hilbert-Schmidt ensembles from the corresponding
results for GUE by using the equivalence of ensembles.
Solving this open problem in the present paper, we use a {\em local form}
of concentration given by the local central limit
theorem for the densities of the squared Hilbert-Schmidt norms.
First we represent the Hilbert-Schmidt measure as  a conditional measure
of the GUE, given the Hilbert-Schmidt norm of the GUE random
matrix. Starting with the disintegration of the GUE according to the
level sets of the Hilbert-Schmidt norm, we arrive at
formula \eqref{bas} which is the crucial ingredient of the proof. Here we have
to extend the scaling parameter to the complex domain. The formal Fourier
inversion applied to this formula gives an 'heuristic proof' of the result. 
For an outline of it see Section 2. To make this sketch rigorous we need asymptotic estimates in the complex
domain for the kernels related to the Hermite functions.  This is done
in Section 3, based on the results in \cite{D99A}, \cite{D99B}, \cite{D99C}. Unfortunately, some of the results we need are contained in these papers not explicitly enough and have to be extracted from the proofs rather than from the statements (see the proof of Lemma \ref{conv}). The main convergence result we need 
is Lemma \ref{conv}. All necessary bounds are summarized in Proposition 
\ref{last1}. Based on these two facts, we complete the proof.  Note that the analytic part of the present paper
may be viewed as a form of the Tauberian theorem.

The authors are indebted to the anonymous referees for careful reading of the manuscript. Moreover, the main theorem in the first version of the paper required 
excluding of a neighborhood of zero (this case has been considered in the note 
\cite{GGL}). Stimulating questions of one of the referees and the editor 
led us to a revision of the proof which removes this restriction. Finally, we would like to thank Justine Swierkot for her assistance in the preparation of the text.

\section{Disintegration, a Fourier transform formula and 
the sketch of the proof}
\label{2}
In this section we discuss  a disintegration  representation of
the GUE in terms of the HSE, and derive a Fourier transform
formula involving these matrix ensembles. We suppress in this
section the "spectral" arguments of the correlation functions 
and related quantities assuming that these arguments vary inside 
the domain described in Section 1.

For every $r>0$ denote by $S_N^r$ the sphere of radius $r$ in
$\HN$ centered at the origin. Let for $s > 0$
$X_N$ be a random matrix in $\HN$ distributed according to
$\m^{\text{GUE},s}_N$. Set $T_N=\tr (X_N^2/s)$ and
$Y_N=NX_N/\sqrt{ \tr X_N^2}  $. The random variable $T_N$ can be
represented as a sum of $N^2$ squares of independent standard
Gaussian random variables, hence it follows the familiar
$\chi^2_{N^2}$ distribution. Moreover, $T_N$ and $Y_N$ are
independent, and $Y_N$ is uniformly distributed on the sphere
$S_N^{N}$ in $\HN$, that is, $Y_N$ is distributed
according to $\nu_N^{1}$ in the notation of the previous section.
 Then $X_N$ can be represented as
\begin{equation} \label{polar1}
X_N=\frac{Y_N}{N}\sqrt{sT_N }
\end{equation}
with $T_N$ and $Y_N$ as above. Let
$\g_{N^2}$ denote the the probability density of $T_N.$ Then it
follows from \eqref{polar1} that
\begin{equation} \label {polar2}
\m_N^{s} =\int_0^{\infty}  \nu_N^{u/N^2} \g_{N^2}(s^{-1}u)
s^{-1} du.
\end{equation}
As a consequence of \eqref{polar2}, the correlation functions of
GUE and HSE for $1 \le n\le N-1$
 satisfy the relation
\begin{equation} \label {discor}
R_{n,N}^{\text{GUE},s} = \int_0^{\infty} R_{n,N}^{\text{HSE},
u/N^2} \g_{N^2}(s^{-1}u) s^{-1} du.
\end{equation}
 Note that
 \begin{equation} \notag
E T_N=N^2 , \ DT_N^2= E(T_N-ET_N)^2 = 2 N^2,
\end{equation}
and, for every $m >0,$
 \begin{equation} \label{chiden}
\g_{m}(u)=\begin{cases}(2^{m/2} \G(m/2))^{-1}x^{(m/2)-1}e^{-x/2}, \ 
\text{if} \, u \, \ge 0, \\
0, \text{if} \, u \, < 0.
\end{cases}
\end{equation}
Set for $s>0$ $\g_{m,s}(\cdot)=s^{-1}\g_{m}(s^{-1}\cdot),$ so
that $\g_{m,1}=\g_{m}$ (note that the same set of densities 
with a different parametrization appears in Lemma \ref{kernelfar} as
$ (f_{a,p})_{a,p>0})$.
Observe now that the probability density of $sT_N$
is given by $ \g_{N^2,s}(\cdot).$ Then \eqref{discor} can be
rewritten as
\begin{equation} \label{dissum}
R_{n,N}^{\text{GUE},s} = \int_0^{\infty} R_{n,N}^{\text{HSE},
u/N^2} \g_{N^2,s}(u) du.
\end{equation}
In particular, we have
\begin{equation} \label {disspec1}
R_{n,N}^{\text{GUE},1/N} = \int_0^{\infty} R_{n,N}^{\text{HSE},
u/N^2} \g_{N^2,1/N}(u) du.
\end{equation}
Our goal is to investigate the limiting behavior of
$R_{n,N}^{\text{HSE}, 1/N} $ when $n$ is fixed, $N \to \infty $,
and the "spectral" arguments of  $R_{n,N}^{\text{HSE}, 1/N} $
vary within an $\e/N-$neighborhood of a point from the main diagonal 
in the cube $(-2,2)^n $. However, we prefer to consider the function $u
\mapsto R_{n,N}^{\text{HSE}, u/N^2 } $ rather than its value
$R_{n,N}^{\text{HSE}, 1/N} $ at $N$. Moreover, we will perform 
the study of the latter function indirectly, first
dealing with the product $u \mapsto R_{n,N}^{\text{HSE}, u/N^2 }
\g_{N^2,1/N}(u).$  After the
change of variable $(u-N)/\sqrt{2}=v $ in \eqref{disspec1} we obtain the 
relation \begin{equation}\label{discen2}
R_{n,N}^{\text{GUE},1/N} 
=\int_{-\infty}^{\infty}q_{N^2}(v) dv, 
\end{equation} where \begin{equation}\label {dens} q_{N^2}(v)
= R_{n,N}^{\text{HSE},1/N +v\sqrt{2}/N^2}\sqrt{2}
\g_{N^2,1/N}(N+v \sqrt{2}). \end{equation}
Notice that  $v \mapsto \sqrt{2}  \ \g_{N^2,1/N}(N+ v \sqrt{2} )
$
 is the probability density of the centered and normalized
random variable $(T_{N^2}-N^2)/\sqrt{2N^2}$, and it tends to the
standard normal density $\varphi: v \mapsto (1/\sqrt{2 \pi})
\exp{(-v^2/2)} $ as $N \to \infty$. Thus, the limit behavior of
the density  $ \g_{N^2,1/N}(N+\cdot \, \sqrt{2})$
 is well understood, and we have to study $q_{N^2}(\cdot).$
For every fixed $u$, due to the relation
$\g_{N^2,s}(\cdot)=s^{-1}\g_{N^2}(s^{-1}\cdot),$ we can analytically extend $\g_{N^2,s}(u)$ to the the domain
$\Re s>0$. Moreover, in the formula \eqref{dissum} the integral
in the right hand side can be analytically continued in $s$ in
accordance with the continuation of $\g_{N^2,s}(u)$ mentioned above. This
leads to the corresponding  continuation of $R_{n,N}^{\text{GUE},s}$
 so that \eqref{dissum}
holds for $s$ from the right half-plane.
Now we will evaluate the Fourier transform of the
(nonprobabilistic) density $q_{N^2}$ keeping in mind that by
\eqref{discen2} and \eqref{dens} it is a nonnegative integrable function.
\begin{lemma}\label{ide}
\begin{equation} \label{bas}
\int_{- \infty}^{\infty}\exp{(ipv)} q_{N^2}(v) dv =
\phi_{N^2}(p) R_{n,N}^{\text{GUE},1/((1-ip\sqrt{2}/N)N)},
\end{equation}
where
\begin{equation} \label{phi}
\phi_{N^2}(p)=\exp{(-ipN/\sqrt{2})} (1-ip\sqrt{2}/N)^{-(N^2/2)}
\end{equation}
is the characteristic function of the random variable
$(T_{N^2}-N^2)/\sqrt{2N^2}$.
\end{lemma}
\begin{proof} Denoting by $C_m$ the normalizing constant in formula (\ref{chiden}), we have
\begin{equation}
\begin{split}
&\int_{- \infty}^{\infty}\exp{(ipv)} q_{N^2}(v) dv \\
=&\int_{- \infty}^{\infty}\exp{(ipv)} R_{n,N}^{\text{HSE},
1/N+v\sqrt{2}/N^2} \sqrt{2} \g_{N^2,1/N}
(N+ v \sqrt{2} )dv \\
=&\int_{- \infty}^{\infty}\exp{\left(\frac{ip(u-N)}{\sqrt{2}}\right)}
R_{n,N}^{\text{HSE}, u/N^2} \g_{N^2,1/N}(u) du\\
=&\exp{\left(\tfrac{-ipN}{\sqrt{2}}\right)}
\int_{- \infty}^{\infty}\exp{(ipu\sqrt{2})}
R_{n,N}^{\text{HSE}, u/N^2}  \g_{N^2,1/N}(u) du \\
=&\exp{\left(\tfrac{-ipN}{\sqrt{2}}\right)} \int_{- \infty}^{\infty}
R_{n,N}^{\text{HSE}, u/N^2}  C_{N^2} N (Nu)^{\tfrac{N^2}{2}-1} 
\exp{\left( -\tfrac{Nu}{2}\Big(1-\tfrac{ip\sqrt{2}}{N}\Big)\right)}du  \\
=&\exp{\left(\tfrac{-ipN}{\sqrt{2}}\right)} \left(1-\tfrac{ip\sqrt{2}}{N}\right)^{-\frac{N^2}{2}+1}\times\\
&\int_{- \infty}^{\infty}R_{n,N}^{\text{HSE}, u/N^2}  
C_{N^2} N \left(Nu\left(1-\tfrac{ip\sqrt{2}}{N}\right)\right)^{\frac{N^2}{2}-1}\!\!\!
\exp{\left( -\tfrac{Nu}{2}\left(1-\tfrac{ip\sqrt{2}}{N}\right)\right)}du \\
=&\exp{\left(\tfrac{-ipN}{\sqrt{2}}\right)} \left(1-\tfrac{ip\sqrt{2}}{N}\right)^{-\frac{N^2}{2}}
\int_{- \infty}^{\infty}R_{n,N}^{\text{HSE}, u/N^2} \g_{N^2,1/((1-ip\sqrt{2}/N)N)}(u) du \\
=& \phi_{N^2}(p) R_{n,N}^{\text{GUE},1/((1-ip\sqrt{2}/N)N)}.
\end{split}
\end{equation}

 \end{proof}
 In the following we will outline our approach. Write
$$ \frac{1}{(Nw(u))^n}q_{N^2}(0)=\frac{1}{2\pi}
\int_{-\infty}^{\infty} \phi_{N^2}(p) \frac{1}{(Nw(u))^n}
R_{n,N}^{\text{GUE},1/((1-ip\sqrt{2}/N)N)}dp. $$
Passing
 to the limit in the integral on the right
hand side (uniformly with respect to the spectral variables),
the right hand side has the same limit as
\[ \frac{1}{(Nw(u))^n}R_{n,N}^{\text{GUE},1/N}\frac{1}{2\pi}
\int_{-\infty}^{\infty}\phi_{N^2}(p)dp\]
 or, in view of the local Central Limit Theorem (CLT) , as
\[\frac{1}{\sqrt{2\pi}(Nw(u))^n}R_{n,N}^{\text{GUE},1/N}. \]
On the other hand, it follows from\eqref{dens} that
\begin{equation}
\frac{1}{(Nw(u))^n} q_{N^2}(0) =\frac{1}{(Nw(u))^n}
R_{n,N}^{\text{HSE}, 1/N} \sqrt{2} \g_{N^2,1/N}(N).
\end{equation}
Again, the local CLT implies that $\sqrt{2} \g_{N^2,1/N}(N) \to
1/ \sqrt{2\pi} $. Therefore,
\[\frac{1}{(Nw(u))^n} R_{n,N}^{\text{HSE}, 1/N}\]
tends to the same limit as
\[\frac{1}{(Nw(u))^n}R_{n,N}^{\text{GUE},1/N},\]
and the conclusion follows.
In the next Section these heuristic arguments will be made rigorous.
\section{Proofs}
\label{3}
For every $\a \in \R,$ through the rest of the paper, we will denote by
$(\cdot)^{\alpha}$ the function
\begin{equation}
(\cdot)^{\alpha}:\C \setminus (-\infty,0] \to \C: z \mapsto
\exp{\a \log z},
\end{equation}
where $\log$ denotes the principal branch of the logarithm. We
use $\sqrt{\cdot}$ as a notation for $(\cdot)^{1/2}$ extended to
$0$ by $\sqrt{0}=0$.

The same symbol (for example, $C, L, M, R, \a,\dots, \d$ with or without indices, $N_0$, and so on) may denote different constants in the bounds obtained in the rest of the paper (even in the same proof). However, we will
indicate explicitly the dependence of such constants on parameters.

 First we are going to state some definitions and basic formulas related to the 
Hermite polynomials and Hermite functions. Then we will formulate some results 
on the asymptotics of the Hermite polynomials and Hermite functions in the
complex plane. The asymptotic behavior of the Hermite
polynomials at the scale we are interested in was first established in 1922 by Plancherel and
Rotach. For a convenient form of these (and much more general) results we refer to 
the monograph  \cite{D99C} and
the papers \cite{D99A},\cite{D99B} (in particular, Appendix B), and \cite{D01}. These results will be 
employed later in this section.

 Let $s > 0$, and let for every $N, N \ge 0,$
 $\tilde{p}_N(\cdot,s)$ be a polynomial of degree $N$ with a
 positive leading coefficient
such that for $(\tilde{p}_N(\cdot,s))_{N \ge 0} $ the relations
$$\int_{\mathbb R}\tilde{p}_M(x,s)\tilde{p}_N(x,s)
\exp{(\! -x^2/(2s))} \, dx =\d_{M,N},\thickspace
 M,N =0,1, \dots,$$ 
 hold.\\
 The (monic) Hermite polynomials
 $(\tilde{H}_N(\cdot,s))_{N=0,1, \dots} $ are defined by the
 expansion
\begin{equation} \label{genfun} 
\exp{(\g x- \g^2s/2)}=\sum_{N=0}^{\infty} \tilde{H}_N(x,s)
\frac{\g^N}{N!}
\end{equation}
and for $M,N =0,1, \dots$ satisfy the relations
$$\int_{\mathbb R} \tilde{H}_M(x,s) \tilde{H}_N(x,s)
\exp{(\! -x^2/(2s))} \, dx =\d_{M,N}
\sqrt{2 \pi} s^{(M+N+1)/2} N!$$ so that
$$\tilde{p}_N(\cdot,s)= \frac{1}{(2\pi)^{1/4} s^{(2N+1)/4}
(N!)^{1/2}}\tilde{H}_N(\cdot,s), \  N=0,1, \dots $$ 
The polynomials
$\tilde{p}_N (\cdot,s)$ satisfy the difference equations
\begin{equation} \label{difeqherm}
 x \, \tilde{p}_N(x,s)= s^{1/2} \, \sqrt{N+1} \,
 \tilde{p}_{N+1}(x,s)+ s^{1/2} \,\sqrt{N}
 \,\tilde{p}_{N-1}(x,s),  \ N=1, 2, \dots
\end{equation}
The Hermite functions
$(\tilde{\varphi}_N(\cdot,s))_{N=0,1,\dots}$ defined by
$$\tilde{\varphi}_N(x,s) \,
= \tilde{p}_N(x,s)\exp{(-x^2/(4s))}, \, N=0,1, \dots, $$
form an orthonormal sequence of functions in
$L_2(\mathbb R,\l)$ where $\l$ is the Lebesgue measure in
$\mathbb R$.
Let us introduce the standardized Hermite polynomials and Hermite functions
(corresponding to the weight $\exp{(-x^2/2)}$) by the relations
\begin{equation} \notag
H_N(\cdot)=\tilde{H}_N(\cdot,1),\,  p_N(\cdot)
= \tilde{p}_N(\cdot,1), \,
\varphi_N(\cdot  )=\tilde{\varphi}_N(\cdot,1)
\end{equation}
so that
\begin{equation} \label{scaling}
\begin{split}
 \tilde{H}_N(\cdot,s)=s^{N/2}&H_N(\cdot \, s^{-1/2}),\, \tilde{p}_N(\cdot,s)= s^{-1/4}p_N(\cdot \, s^{-1/2}),\\ 
&\tilde{\varphi}_N(\cdot,s)=s^{-1/4} \varphi_N(\cdot \, s^{-1/2}).
\end{split}
\end{equation}
The above relations for the Hermite polynomials
imply that the Hermite functions satisfy,
for every $k \ge 1,$ the following system of differential
equations:
  \begin{equation} \label{diffeq}
\begin{split}
&\varphi_k'(x)=-\frac{x}{2} \varphi_k(x)
+ \sqrt{k} \varphi_{k-1}(x), \\
&\varphi_{k-1}'(x)=-\sqrt{k} \varphi_k(x) +\frac{x}{2}
\varphi_{k-1}(x).
\end{split}
\end{equation}
The reproducing kernel $\tilde{K}_N(\cdot,\cdot,s)$ of the
orthogonal projection in $L_2(\R, \l)$ onto the linear span of
$\tilde{\varphi}_0(\cdot,s),
\dots,\tilde{\varphi}_{N-1}(\cdot,s)$ is given by
\begin{equation} \label{kerexp}\tilde{K}_N(x,y,s)=\sum_{k=0}^{N-1} \tilde{\varphi}_k(x,s)\tilde{\varphi}_k(y,s).
\end{equation}
Note that the correlation functions of the GUE can be expressed in terms of the latter kernel as 
\begin{equation} \label{gdet}
R_{n,N}^{\text{GUE},s}(x_1,\dots,x_n)=\det  \biggl( \tilde{K}_N(x_i,x_j,s)\biggr)_{\mspace{-4.0 mu} i,j=1}^{\mspace{-4.0 mu}n}  \mspace{-16.0 mu}.
\end{equation}
(see \cite{M}, (6.2.7)).

Further, setting 
$$ K_N(x,y)= \tilde{K}_N(x,y,1),$$ 
we obtain
\begin{equation} \label{sckernel}
\tilde{K}_N(x,y,s)= s^{-1/2} K_N(xs^{-1/2},ys^{-1/2}).
\end{equation}
The Christoffel - Darboux formula for the kernels $\tilde{K}_N$ and $K_N$ 
reads
\begin{equation}\label{cdsc}
\tilde{K}_N(x,y,s)= \frac{\tilde{\varphi}_N(x,s)\tilde{\varphi}_{N-1}(y,s) - \tilde{\varphi}_{N-1}(x,s) 
\tilde{\varphi}_N(y,s)}{x-y}
\end{equation} 
and 
\begin{equation}\label{cdnrm}
K_N(x,y)=\frac{\varphi_N(x)\varphi_{N-1}(y) - \varphi_{N-1}(x)\varphi_N(y)}{x-y}.
\end{equation} 
Note that the kernels $\tilde{K}_N, K_N$ enjoy the property
\begin{equation} \label{symm}
 \tilde{K}_N(-z_1,-z_2,s)=\tilde{K}_N(z_1,z_2,s), \, K_N(-z_1,-z_2)=K_N(z_1,z_2).
\end{equation} 
As to $K_N,$ this is a consequence of \eqref{cdnrm} and the fact that the functions 
$(\varphi_N(\cdot))_{N \ge 0}$ are even or odd depending on $N$ being even or odd;  for  $\tilde{K}_N$ it follows from \eqref{sckernel}.
The following integral representation of the reproducing kernel
is a version of formula (4.56) in \cite{Forr}:
\begin{equation} \label{Forr}
\begin{split}
& K_N(x,y) = \\
&\sqrt{\frac{N}{2}}
\int_0^{\infty}(\varphi_N(x+\ta)\varphi_{N-1}(y+\ta)
+\varphi_{N-1}(x+\ta)\varphi_N(y+\ta))d\ta
\end{split}
\end{equation}
so that
\begin{equation} \label{scForr}
\begin{split}
&\tilde{K}_N(x,y,s)= \\
&\sqrt{\frac{N}{2s}}\int_0^{\infty}\!\!\!
\left(\varphi_N\bigl(\tfrac{x}{\sqrt{s}}+\tau \bigr)\varphi_{N-1}\bigl(\tfrac{y}{\sqrt{s}}+\ta\bigr) 
+\varphi_{N-1}\bigl(\tfrac{x}{\sqrt{s}}+\ta\bigr)\varphi_N\bigl(\tfrac{y}{\sqrt{s}}+\ta\bigr)\right)d\ta.
\end{split}
\end{equation}
The relations \eqref{scaling}, \eqref{sckernel}, and our extension of
the function $(\cdot)^{\a}$ to $\mathbb C \setminus (- \infty,
0)$ allow us to continue $ \tilde{H}_N(x,\cdot)$, $ \tilde{p}_{N}(x,\cdot),$ $
\tilde{\varphi}_{N}(x,\cdot)$
 and $\tilde{K}_N(x,y,\cdot)$ to this domain analytically in the parameter $s.$ The relations \eqref{genfun}\,-\,\eqref{symm}
 remain valid under these continuations. 
 So does (\ref{scForr}) whenever the integrals there are
 well-defined.
 
Denote by $w$ the analytic continuation of the standard 
Wigner density \eqref{wigner} to the domain $\C \setminus
\bigl((-\infty,2] \cup [2, \infty)\bigr).$
Let us define for every $\a \in (0, 2)$ and  $\b >0$ the sets
$S_{\a,\b}$ and $\overline{S_{\a,\b}}$ by the relations
\begin{equation} \notag
S_{\a,\b}= \{z\in \C:|\Re z| < 2-\a, |\Im u| < \b \}
\end{equation}
and 
\begin{equation} \notag
\overline{S_{\a,\b}}= \{z\in \C:|\Re z| \le 2-\a, |\Im u| \le \b \}.
\end{equation}
Set for $ H\in\R$ \,$d(H)=\sqrt{1+iH}$ and observe that
\begin{equation} \label{scalmod}
\hspace{-1.55cm}|d(H)| =(1+H^2)^{1/4},
\end{equation}
\begin{equation} \label{scalre}
\hspace{-0.95cm} \Re d(H)= \sqrt{\frac{\sqrt{1+H^2}+1}{2}},
\end{equation}
\begin{equation} \label{scalim}
\Im d(H)= \sgn H \,\sqrt{\frac{\sqrt{1+H^2}-1}{2}} ,
\end{equation}
and
\begin{equation} \label{scalhyp}
\hspace{-1.53cm}\bigl(\Re d(H)\bigr)^2- \bigl(\Im d(H)\bigr)^2=1.
\end{equation}
It is clear from \eqref{scalhyp} that
\begin{equation} \label{cone}
|\Re(u d(H))|\ge |\Im( ud(H))|
\end{equation}
for every $u \in \R.$
Note that for every $b \ge 0$ the equation
$$ |\Im d(H)|= b $$
has a unique nonnegative solution
\begin{equation} \label{Hbound}
 H_b=\sqrt{(1+2b^2)^2-1}.
\end{equation}
\begin{lemma} \label{switchone} Let $\a,$ $\b,$ and $u$ be real numbers 
such that $\a \in (0,1),$ $0<\b <\a,$ and $u \in (-2+2\a,2-2\a)$. Then the following assertions hold true: 
\begin{enumerate} 
\item For every  $b \in [-\infty, \infty]$ the relations $|H|< |H_b|$ and $|\Im \bigl(d(H)\bigr)| < b$ are equivalent (we set here $H_{\infty}=\infty,$ 
$\Im(d\bigl(\infty)\bigr)=\infty,$ and $\Im\bigl(d(-\infty)\bigr)=-\infty$); the  relations $|H|\le |H_b|$ and $|\Im \bigl(d(H)\bigr)| \le b$ are also equivalent.
\item For every real $H$ we have $ud(H) \in \overline{S_{\a,\b}}$ whenever   $|\Im(ud(H))| \le \b$;  
\item If for a certain real $H$ the relation $|H|\le |H_{\b/2}|$ (or, equivalently, $|\Im \bigl(d(H)\bigr)| \le \b/2$) 
holds, it follows that $ud(H) \in \overline{S_{\a,\b}}.$  
\end{enumerate}
\end{lemma}
\begin{proof}
It is clear from \eqref{Hbound} and the definition of $H_b$ that for $b \in 
[0, \infty]$ and $H \in [0,\infty]$ 
the functions $b \mapsto H_{b}$ and $H \mapsto |\Im d(H)|$ are strictly increasing mutually inverse 
mappings. This shows equivalence of the conditions in point 1
of the lemma for the functions in question restricted to $[0,\infty]$. The same is true 
for $b, H \in [-\infty, \infty]$ since both functions are even. This proves assertion 1.\\ 
 For every $ H$  such that \begin{equation} \label{ineq2}
\bigl|\Im \bigl(ud(H)\bigr)\bigr| \le \b
\end{equation} 
we have
\begin{equation} \label{ineq1}
\bigl|\Re \bigl(ud(H)\bigr) \bigr|  < 2-\a
\end{equation}
since it follows from the relation \eqref{scalhyp} and the assumptions $\b\le \sqrt{\a}$, $u \in (-2+2\a,2-2\a)$ that 
 \begin{equation} \notag
\begin{split}
\bigl|\Re \bigl(ud(H)\bigr) \bigr|
= \sqrt{u^2+\bigl(\Im u\bigl(d(H)\bigr)\bigr)^2} 
& \le \sqrt{(2-2\a)^2+ \b^2}\\ 
 &<\sqrt{(2-2\a)^2+ \a^2}
<2-\a.
\end{split}
\end{equation}
In view of the definition of $\overline{S_{\a,\b}}$, the second assertion is proved. 
The third assertion is an immediate consequence of the second one since $|u| \le 2.$ 

\end{proof}
\begin{corollary} \label{switchone1}Let $\a,$ $\b,$ and $u$ satisfy the assumptions of Lemma \ref{switchone}, and assume that the relation $ud(H) \notin \overline{S_{\a,\b}}$ holds for some real $H.$  Then we have
$|\Im(u d(H))|> \b.$ Moreover, the relation $|\Im(ud(H'))| > \b$ holds for every real $H'$ 
such that $|H'|\ge |H|.$  
\end{corollary}
\begin{proof}
The corollary follows immediately from the second 
and the first assertions of Lemma \ref{switchone} combined with the fact that, by \eqref{scalim}, 
$|\Im d(H)|$ is a strictly 
increasing function of $|H|.$ 

\end{proof}
In the next Lemma the main result on convergence is formulated.
\begin{lemma} \label{conv} For every number $\a \in (0,1)$ there exists a positive number $\overline{H}=\overline{H}(\a)$ 
such that for every $A>0$ the relation
\begin{equation} \label{estim1}
\begin{split}
&\frac{1}{Nw(u)}\tilde{K}_N
\biggl(\biggl(u+\frac{t_1}{Nw(u)}\biggr)d(H),
\biggl(u+\frac{t_2}{Nw(u)}\biggr)d(H), N^{-1}\biggr)\\
&\underset{N\to \infty}{\longrightarrow} \frac{\sin
\bigl(\pi(t_1-t_2)d(H)w(ud(H))/w(u)\bigr)}{ \pi(t_1-t_2)d(H)}
\end{split}
\end{equation}
holds uniformly for all real numbers $u,$ $t_1, t_2,$ and $H$ 
provided that $u \in [-2+2\a,2-2\a],$ $|t_1| \le A,$ $|t_2| \le A,$ and  
$|H| \le \overline{H}.$
\end{lemma}
\begin{proof}
Note that for every $\a\in(0,1)$ there exists
$\b'=\b'(\a)>0$ such that for every $A > 0$ the relation
\begin{equation} \label{sinlim}
\frac{1}{Nw(v)}\tilde{K}_N\biggl (
v+\frac{z_1}{Nw(v)},v+\frac{z_2}{Nw(v)}, \frac{1}{N}
\biggr)\underset{N\to \infty}{\to}
 \frac{\sin \pi(z_1-z_2)}{ \pi(z_1-z_2)}
\end{equation}
holds uniformly for all complex numbers $v \in
\overline{S_{\a,\b'(\a)}},$ $z_1,$ and $z_2 $ such that  $| z_1| \le A, |z_2|
\le A.$\\
 This assertion (and more general ones involving
some class of weights) is impli\-citly contained in the papers \cite{D99A} and
\cite{D99B} (see also the monograph \cite{D99C}). Actually, Lemma 6.1 in
\cite{D99A} establishes the desired result for real $u,$ $
z_1$ and $z_2.$ However, the same reasoning applies to the complex $u,$ $
z_1,$ and $z_2$
satisfying the assumptions just made provided that
$\b'$ is chosen to be sufficiently small in accordance with $\a$. The key 
ingredient of the proof 
in these references, the boundedness property of a
certain derivative, is established in  \cite{D99A} and
 \cite{D99B} for some complex neighborhood of an arbitrary real point
 $u \in [-2+2\a, 2-2\a]$ (see relation (4.122) in \cite{D99A})
 which allows to bound
this derivative in a rectangular strip $ \overline{S_{\a,\b'(\a)}}.$ Lowering $\b'$ if
necessary,  we can substitute $\b'(\a)$ by some $\b(\a) \le
\a/4.$ Note that the function $w(\cdot)$ has no zeroes in $\C \setminus
\bigl((-\infty,2] \cup [2, \infty)\bigr)\supset
\overline{S_{\a,\b(\a)}}$.
Picking a point $v' \in \overline{S_{\a,\b(\a)}}$ and setting in \eqref{sinlim} $z_i=z'_i \, w(v)/w(v') \,  \, \,
(i=1,2),$ we obtain
\begin{equation} \label{sinlim1}
\frac{1}{Nw(v')}\tilde{K}_N\biggl (
v+\frac{z'_1}{Nw(v')},v+\frac{z'_2}{Nw(v')}, \frac{1}{N}
\biggr)\underset{N\to \infty}{\to}
\frac{\sin \biggl(\frac{\pi(z'_1-z'_2)w(v)}{w(v')}\biggr)}
{ \pi(z'_1-z'_2)}.
\end{equation}
Note that  for $v \in \overline{S}_{\a,\b(\a)}$ we have $ C^{-1}\le|w(v)|\le C$
with some $C> 1.$ This fact and \eqref{sinlim} imply that for every 
$A > 0$ \eqref{sinlim1}
holds uniformly in $v, v' \in \overline{S_{\a,\b(\a)}}$  and $| z'_1| \le A,
|z'_2| \le A.$

According to assertion 3 of 
Lemma \ref{switchone}, for every $u \in [-2+
2\a, 2 - 2\a ]$  and every $H \in [-H_{\b(\a)/2},
H_{\b(\a)/2}] $ we have 
$ud(H) \in \overline{S_{\a,\b(\a)}}.$
 More, it is clear from \eqref{scalmod} that
for $H \in [-H_{\b(\a)/2},
H_{\b(\a)/2}] $ we have ${C^{'}}^{-1}<|d(H)|<C'$ with some positive constant $C'>1$ depending on $\a.$
Then, setting $v'=u \in [-2+2\a,2-2\a],$ $v= u d(H) \in  \overline{S_{\a,
\b(\a)}},$ $z'_i= t_i d(H)$ $(i=1,2),$ we obtain  \eqref{estim1}.

\end{proof} 
\begin{remark} \label{conv1} 
It is sometimes convenient to use the following reformulation of \eqref{sinlim} and the related 
assumptions :\\
{\em For every number $\a \in (0,1)$ there exists a real numbers $\b=\b(\a)\le \a/4$  
such that for every $A>0$ the relation
\begin{equation} 
\notag
\biggl|\frac{1}{N}\tilde{K}_N
\bigl(v_1,v_2, N^{-1}\bigr)\\
 - \frac{\sin
\bigl(N \pi(v_1-v_2)w(v)\bigr)}{ N\pi (v_1-v_2)}\biggr|\underset{N\to \infty}{\longrightarrow} 0
\end{equation}
holds uniformly for all complex numbers $v,$ $v_1, v_2$ satisfying the 
constraints  
 $v \in \overline{S_{\a,\b(\a)}},$ 
$|v_1-v| \le A/N,$ $|v_2-v| \le A/N$. }\\

\end{remark}
 We will now derive upper bounds for the modulus of the kernel
\begin{equation} \notag
\frac{1}{N}\tilde{K}_N
\bigl(\cdot,
\cdot, N^{-1}\bigr)
\end{equation}
 valid for various values of its arguments. We will use asymptotics of two types for the Hermite functions: one of them is valid in a narrow strip around the interval $[-2+\a, 2-\a];$ another one applies in the 
complement to some neighborhood of $[-2,+2].$ With this asymptotics, after obtaining certain intermediate bounds for the quantity $|N^{-1} \tilde{K}_N
\bigl(u+ t_1/N)d(H),u+ t_2/N)d(H), N^{-1}\bigl)|$ , we will finally determine for this quantity an upper estimate which does not depend on $H \in \R.$ This bound and analogous bounds for the correlation functions will be derived in Proposition \ref{last1}.\\
 In the next two lemmas we present upper estimates for the Hermite functions. These estimates are immediate consequences of the Plancherel-Rotakh-type asymptotics. First we consider a bound valid in the strip $\overline{S}_{\a,\b(\a)}$ where $\b(\a)$ was defined in the proof of Lemma \ref{conv}. 
\begin{lemma} \label{hermfunnear}For any $\a \in (0,1)$ there exist a constant $C(\a) > 0$ such that we have
\begin{equation}
\label{hermfunNnear} \bigl|\varphi_N(z \sqrt{N})\bigl| \le
\bigl(C(\a)\bigl)^N
\end{equation}
and
 \begin{equation} \label{hermfunN-1near} 
\bigl|\varphi_{N-1}\bigl(z \sqrt{N}\bigr)\bigr|
\le \bigl(C(\a)\bigl)^{N-1}
\end{equation}
for every  $N \in \N $ and  $z \in \overline{S}_{\a,\b(\a)}.$
\end{lemma}
\begin{proof}
First we recall a known result about the Plancherel-Rotach
asymptotics for the Hermite functions in a rectangular strip  
(Theorem 2.2, part (ii), in \cite{D99B}; the particular case of the Hermite functions for the measure $\exp{(-x^2)} dx$ is considered there in Appendix B). This strip  in the complex plane 
contains, up to some neighborhoods of the ends, the interval where the zeroes of the Hermite functions asymptotically concentrate.
Set 
\begin{equation} \label{}
\notag
\psi: \C \setminus((-\infty,-1] \cup[1,\infty)) \to \C: z\mapsto \frac{1}{2 \pi}(1-z)^{1/2}(1+z)^{1/2}.
\end{equation}
The function $\arcsin$ is defined  as the inverse function of 
\[ \notag 
\sin: \biggl \{ z \in \C: |\Re (z)| < \frac{\pi}{2} \biggr \} \to 
\C \setminus((-\infty,-1] \cup[1,\infty)).
\]
In our notation, the statement in \cite{D99B} reads: 
\nolinebreak
{\em for every $\a \in (0,1) $ uniformly in  $z \in$ $\{z' \in \C: |\Re(z')| \le 1-\a,|\Im(z')| \le \a \}$
we have
\begin{equation}
\label{asymp}
\begin{split}
\tilde{\varphi}_N(\sqrt{2N}z,1/2) &\\
=\sqrt{\frac{2}{\pi \sqrt{2N}}} (1&-z)^{-1/4}(1+z)^{-1/4} \\
\times \biggl\{ \cos \biggl( &N \pi  \int_1^z \psi (y)dy + \frac{1}{2} \arcsin z \biggr) \biggl( 1+O\biggl(\frac{1}{N}\biggr)\biggr) \\
&+\sin \biggl(  N \pi  \int_1^z \psi (y)dy -\frac{1}{2} \arcsin z  \biggr)
O\biggl(\frac{1}{N}\biggr)\biggr\}.
\end{split}
\end{equation}
}
The analytic functions $ z \mapsto \int_1^z \psi (y)dy$ and $\arcsin $ are bounded on the set $ \{z' \in \C: |\Re(z')| \le 1-\a,|\Im(z')| \le \a \}.$ Representing $\cos$ and $\sin$ as the sum and the difference of the exponentials, it follows from \eqref{asymp}\,(omitting the multiplier $N^{-1/4}$) that the modulus of the left hand side of \eqref{asymp} is bounded above by the $N-$th power of 
a positive number depending on $\a.$  For every $\a \in (0,2),$ due to the relations
$$\varphi_N ( z \sqrt{N}) 
= 2^{-1/4} \tilde{\varphi}_N\biggl(\frac{z}{\sqrt{2}}\sqrt{N},1/2\biggl)
= 2^{-1/4} \tilde{\varphi}_N\biggl(\frac{z}{2}\sqrt{2N},1/2\biggl), $$ 
a similar bound by the $N-$th power of some number depending on $\a$ holds for $\varphi_N( z \sqrt{N})$ uniformly for $z \in \overline{S}_{\a,\a}.$
For 
$$\varphi_{N-1} ( z \sqrt{N})=\varphi_{N-1} \bigl( z\sqrt{(N-1)}\sqrt{N/(N-1)} \bigr)$$
we also have a bound of such type in $\overline{S}_{\a,\a/4} \subset \overline{S}_{\a/2,\a/2}.$  Indeed, set $\rho_N=\sqrt{N/(N-1)}$, 
and let $N_0=N_0(\a)$ be defined by the equation $\rho_{N_0} = (2-\a/2)/(2-\a).$
Then we have $\rho_N \le 2$ and $\rho_N \overline{S}_{\a,\a/4} \subset \overline{S}_{\a/2,\a/2} $
whenever $N \in \N$ satisfies $N \ge N_0.$  This means that for  $z \in  \overline{S}_{\a,\a/4}$ both $z \in  \overline{S}_{\a/2,\a/2}$ and $z\sqrt{N/(N-1)}\in  \overline{S}_{\a/2,\a/2}$ hold for every $N \ge N_0.$ Therefore, the conclusion of the Lemma  is valid for such values of $N.$ Choosing a larger constant $C(\a)$ if necessary and recalling that $\b(\a) \le \a/4$, the proof is completed.

\end{proof} 
We will obtain now some estimates valid the Hermite functions in the domain $|\Im z| \ge \d >0.$ 
\begin{lemma} \label{hermfunfar}
For every $\d >0$ there exist constants \,$ C(\d)$ and
$M(\d)$ such that the inequalities
\begin{equation}
\label{funcstripN} \bigl|\varphi_N(z \sqrt{N})\bigl| \le
C(\d)N^{-1/4}M^N(\d)\,|z|^N
\end{equation}
and
 \begin{equation}
\label{funstripN-1} \bigl|\varphi_{N-1}\bigl(z \sqrt{N}\bigr)\bigr|
\le C(\d)(N)^{-1/4}M^{N-1}(\d)\,|z|^{N-1}
\end{equation}
hold for every $N \in \N$ and every $z$ satisfying $|\Re z|
\ge |\Im z| \ge \d.$
\end{lemma}
\begin{proof}
Again we start with stating some known result about the Plan\-che\-rel-Rotach
asymptotics for the Hermite polynomials in the complex plane 
(Theorem 7.185 in \cite{D99C}, see also \cite{VAG89} and the
references therein). In our notation, uniformly in  $z$
from every compact set contained in $(\C \cup\{ \infty \})
\setminus [-\sqrt{2},\sqrt{2}],$ we have
\begin{equation}
\label{asymp1}
\tilde{H}_N(z,(2N)^{-1}) =\frac{U(z)}{2}\frac
{\exp{(N(z-\sqrt{z^2-2})^2}/4)}{(z-\sqrt{z^2-2})^N}
\biggl(1+O\biggl(\frac{1}{N}\biggr)\biggr),
\end{equation}
where 
\begin{equation} 
\notag
U(z):=\biggl(\frac{z-2}{z+2})^{1/4}
+(\frac{z+2}{z - 2}\biggr)^{1/4}.
\end{equation}
From this relation we obtain for the monic orthogonal Hermite polynomials $(H_N)_{N\ge 0}$ with respect
to  the weight $\exp{(-x^2/2)}$ that 
\begin{equation}
\label{asymp2}
\begin{split}
H_N(z\sqrt{N})&=(2N)^{N/2}\tilde{H}_N(z/\sqrt{2},(2N)^{-1}) \\
&=2^{N-1}N^{N/2}U(z)\frac
{\exp{(N(z-\sqrt{z^2-4})^2}/8)}{(z-\sqrt{z^2-4})^N}
 \biggl(1+O\biggl(\frac{1}{N}\biggr)\biggr)
\end{split}
\end{equation}
uniformly for $z$ from every compact set contained in $(\C \cup\{ \infty \})
\setminus [-2,2].$ Passing to the normalized polynomials we see that
\begin{equation}
\begin{split}
p_N(z \sqrt{N})
&=\frac{1}{(2 \pi)^{1/4}(N!)^{1/2}} H_N(z \sqrt{N}) \\
&=\frac{2^N N^{N/2}U(z)}{2(2 \pi)^{1/4}(N!)^{1/2}}
\frac
{\exp{(N(z-\sqrt{z^2-4})^2}/8)}{(z-\sqrt{z^2-4})^N} \biggl(
1+O\biggl(\frac{1}{N}\biggr)\biggr),
\end{split}
\end{equation}
and, in view of the Stirling formula $ N!= \sqrt{2 \pi}
N^{N+\frac{1}{2}}e^{-N} (1+O(1/N)),$ we obtain
\begin{equation}
\label{asymp3}
\begin{split}
p_N(z \sqrt{N}) = 
\frac {U(z)\, \exp{(N( \frac{1}{2}+ \log 2+
(z-\sqrt{z^2-4})^2}/8))} {2(2 \pi)^{1/2}(N)^{1/4}\, (z-\sqrt{z^2-4})^N} \biggl(
1+O\biggl(\frac{1}{N}\biggr)\biggr)
\end{split}
\end{equation}
so that we have
\begin{equation}
\label{asymptout}
\begin{split}
p_N(z \sqrt{N})
 = \frac {U(z)\,\exp{( \frac{N}{2}(1+ ((z-\sqrt{z^2-4})/2)^2}))}
{2(2 \pi)^{1/2}(N)^{1/4}\bigl((z-\sqrt{z^2-4})/2\bigl)^N} \biggl(
1+O\biggl(\frac{1}{N}\biggr)\biggr)
\end{split}
\end{equation}
uniformly in $z$ from any compact subset of $(\C \cup\{ \infty
\}) \setminus [-2,2].$\\
 For every complex number $z \in \C \setminus [-2,2]$ let $x=x(z)$ be the root
of the equation $x+x^{-1}=z$ satisfying $|x| < 1.$ Then we have
$$ x=x(z)= \bigl(z-\sqrt{z^2-4}\bigr)/2 $$
where $z \mapsto \kappa(z)=\sqrt{z^2-4}$ is, by definition, a
univalent analytic function  in $\C \setminus [-2,2]$ satisfying
the relations $\kappa^2(z)=z^2-4$ and $\kappa(t) > 0$ for $t>2.$
For every $r \in (0,1)$ the equation $ |x(z)|= r$ and the inequality $ |x(z)|\le r$  define the ellipse 
\begin{equation} \label{ellipse} (\Re z)^2/(r^{-1}+r)^2+ (\Im z)^2/(r^{-1}-r)^2=1
\end{equation}
with focal points $-2,2$ and the closed exterior domain $E_r$ of this ellipse, respectively.  Thus for $z \in E_r$ we have
\begin{equation} \notag
|x(z)|= \bigl|\bigl(z-\sqrt{z^2-4}\bigr)/2\bigr|\le r < 1,
\end{equation}
and
\begin{equation} \label{numer}
 \Re \bigl((z-\sqrt{z^2-4})/2\bigr)^2 \le r^2 <1.
\end{equation}
Note that for every $r \in (0,1)$
\begin{equation} \notag
\min_{|x(z)|=r}|z|= r^{-1}-r
\end{equation}
which  for $ z \in E_r$ implies
\begin{equation} \label{denom}
\begin{split}
\frac{2}{\bigl|z-\sqrt{z^2-4}\bigr|}&
=\frac{1}{|x(z)|}=|z-x(z)| \le |z|+1 \le (1+|z|^{-1})|z| \\
&\le \bigg(1+\frac{1}{ r^{-1}-r}\biggr)|z|.
\end{split}
\end{equation}
Further, for $z \in E_r$ $(0<r<1) $ we can write
\begin{equation} \notag
\begin{split}
|z-2|& 
=\bigl|\bigl(x(z)-1\bigr)+\bigl(x^{-1}(z)-1\bigr)\bigr|
= 2|x(z)-1||x^{-1}(z)-1| \\
&=2|x(z)|^{-1}|1-x(z)|^2  \ge 2r^{-1}|1-x(z)|^2 \ge 2r^{-1}|1-r|^2 \\
&=2(r^{-1/2}-r^{1/2})^2,
\end{split}
\end{equation}
and analogously
\[ \notag |z+2| \ge  2(r^{-1/2}-r^{1/2})^2. \]
This implies that the function
\begin{equation} \notag
|U(z)|=\biggl|\biggl(\frac{z-2}{z+2}\biggr)^{1/4}
+\biggl(\frac{z+2}{z-2}\biggr)^{1/4}\biggr|
\end{equation}
is bounded above on $E_r$ by a constant depending on $r \in
(0,1).$ Applying this estimate along with \eqref{numer} and \eqref{denom} 
to the relation \eqref{asymptout} we
arrive, for a fixed $ r \in (0,1)$ and every $ z \in E_r,$ at
\begin{equation} \label{polyN}
 \bigl|p_N(z \sqrt{N})\bigl| \le C_1(r)N^{-1/4}L^N(r)|z|^N.
\end{equation}
 Now we set $z_n=z \sqrt{N/(N-1)}$ and note that, by convexity
 of $\C\setminus E_r,$ $ z_N \in E_r$ if $z \in E_r.$
Hence,  applying  \eqref{polyN} to $p_{N-1},$ we see that for
every $z \in E_r$
\begin{equation} \notag
\begin{split}
\bigl|p_{N-1}\bigl(z \sqrt{N}\bigr)\bigr|&
=\bigl| p_{N-1}\bigl(z_N \sqrt{N-1}\bigr)\bigr| \, \le  \, C_1(r)(N-1)^{-1/4}L^{N-1}(r)|z_N|^{N-1} \\
 &= C_1(r)\bigl(1+1/(N-1)\bigr)^{(2N-3)/4}N^{-1/4}
 L^{N-1}(r)|z|^{N-1}.
\end{split}
\end{equation}
Since $\bigl(1+1/(N-1)\bigr)^{(2N-3)/4} \to e^{1/2}$ as $N
 \to \infty,$ we conclude that
\begin{equation}
\label{polyN-1} \bigl|p_{N-1}\bigl(z \sqrt{N}\bigr)\bigr|\le
C_2(r)(N)^{-1/4}L^{N-1}(r)|z|^{N-1}.
\end{equation}
Observe that  for every $r \in (0,1)$ the inequality
\[ 
 \Im z \ge r^{-1} - r 
 \]
implies, by equation \eqref{ellipse}, that
$z \in E_r.$
 Finally set
\[
r =-\frac{\d}{2}+\sqrt{\frac{\d^2}{4}+1}
\]
so that $r$ is the solution of the equation
$r^{-1}-r=\d$
satisfying $0<r<1$.
We can now conclude from \eqref{polyN} and \eqref{polyN-1} that
for every $\d >0$ with some constants $ C(\d)$ and
$M(\d)$ the inequalities
\begin{equation}
\label{polystripNfar} \bigl|p_N(z \sqrt{N})\bigl| \le
C(\d)N^{-1/4}M^N(\d)|z|^N
\end{equation}
and
 \begin{equation} \label{polystripN-1far} 
\bigl|p_{N-1}\bigl(z \sqrt{N}\bigr)\bigr|
\le C(\d)(N)^{-1/4}M^{N-1}(\d)|z|^{N-1}
\end{equation}
hold for every natural $N$ and every $z$ with $|\Im z| \ge \d.$
We pass now from Hermite polynomials to Hermite functions.
Note that for every $z$ satisfying $|\Re z|\ge|\Im z|$ we have
\begin{equation} \label{boexp}
|\exp{(-z^2/4)}| \le \exp{\bigl((-\Re^2 z+ \Im^2z)/4\bigr)} \le 1
\end{equation}
and, therefore,
\begin{equation} \label{bovar}
| \varphi_k(z)| \le |p_k(z)|
\end{equation}
for every $k \in \N .$
Applying \eqref{polystripNfar} and \eqref{polystripN-1far}, the Lemma follows.
\end{proof}
\begin{corollary} \label{pathfar}
For every real numbers $\d >0$ and $A>0$ there exist such constants \,$ C(\d,A)$ and
$M(\d)$ that the inequalities
\begin{equation}
\label{funcstrippathN} \bigl|\varphi_N\bigl((u+t/N)\,d(H) \sqrt{N}\bigr)\bigl| \le
C(\d,A)M^N(\d)|H|^N
\end{equation}
and
 \begin{equation}
\label{funcstrippathN-1} 
\bigl|\varphi_{N-1}\bigl((u+t/N)\,d(H) \sqrt{N}\bigr)\bigl| \le
C(\d,A)M^{N-1}(\d)|H|^{N-1}
\end{equation}
hold for every numbers $u \in [-2,2],$ $t \in [-A,A],$ $N \in \N$ and $H \in \R$ such that 
$\bigl|\Im\bigl({(u+t/N)d(H)}\bigr)\bigr|\ge \d.$
\end{corollary}
\begin{proof}
Set $ v=(u+t/N)\,d(H).$ 
Since
\begin{equation} \notag 
|\Im(d(H))|=\sqrt{\bigl(\sqrt{1+H^2}-1\bigr)/2} \le |H|/2, 
\end{equation}
we have that
\begin{equation} \notag
|\Im(v)| \le (2+(A/N))|\Im({d(H)})|  \le
\bigl(1+ (A/2N)\bigr)|H|
\end{equation}
and
\begin{equation}  \notag
|\Im(v)|^N  \le
 \bigl(1+(A/2N)\bigr)^N |H|^N\le  C(A) |H|^N,
\end{equation}
since
\begin{equation} \notag
\bigl(1+(A/2N)\bigr)^N  \underset{N \to \infty}{\to}\exp{(A/2)}.
\end{equation}
Thus, we proved that 
\begin{equation} \label{expbound}
|\Im(v)|^N  \le  C(A) |H|^N.
\end{equation}
In view of \eqref{cone},  $|\Re v| \ge |\Im v|,$ and we can apply
Lemma \ref{hermfunfar}. Omitting $N^{-1/4}$ in the bound, the proof is 
completed. 

\end{proof}
Now we will derive upper bounds for $ |N^{-1}\tilde{K}_N(\cdot, \cdot, N^{-1})|$  valid for various combinations of the values of its arguments. 
\begin{lemma}\label{kernelnear} For every real $\a \in (0,1)$ and $A>0$ there exists
 a constant $R_1(\a,A)>0$ such that for every real $ u \in [-2+\a,2-\a],$ $t_1, t_2 \in [-A,A],\, H,$ and every  $N \in \N$ we have
\begin{equation} \label{bound}
|N^{-1}\tilde{K}_N \bigl((u+t_1/N)d(H),(u+t_2/N)d(H), N^{-1})|  
\le \bigl(R_1(\a,A)\bigr)^{2N}
\end{equation}
whenever 
\begin{equation} \label{inside}
|\Im(u+t_i/N) d(H)| \le \b(\a/4)
\end{equation}
holds for $i=1,2.$
\end{lemma}
\begin{proof}
Set 
$$u_i=u+t_i/N,\quad  v_i = (u+t_i/N)d(H)\quad \text{for} \, \, i=1,2 $$ 
and 
$$u= (u_1+u_2)/2, v = (v_1+v_2)/2.$$

Let the number $N_0=N_0(\a, A)$ be large enough such that the inequality
 \begin{equation} \label{nefirst}
A/N \le \a/2  
\end{equation}
holds for every $N \ge N_0.$ 
Observe that by \eqref{nefirst}
$ u_i \in [-2+\a/2,2-\a/2], i =1,2, $ 
for $N \ge N_0.$
Then it follows from \eqref{inside} and Corollary \ref{switchone1} 
that $ v_1, v_2 \in \overline{S_{\a/4,\b(\a/4)}}.$ 
Therefore, $v \in \overline{S_{\a/4,\b(\a/4)}},$ too. \\
First assume now that we have 
\begin{equation} \label{narrow}
 |v_1 -v_2| \le 1/N.
\end{equation}
Then Remark \ref{conv1} (with $A=1$) applies, and it is clear that, 
for $v_1,v_2,$ $v =$ \linebreak 
$(v_1+v_2)/2 \in \overline{S_{\a/4,\b(\a/4)}}$ 
 with $v_1,v_2$ satisfying \eqref{narrow}, the function \linebreak
$|\sin
\bigl(N \pi(v_1-v_2)w(v)\bigr)/\bigl( N\pi (v_1-v_2)\bigr)|$
is bounded above by a constant depending on $\a$ and $A.$ Therefore, by Remark \ref{conv1}, this applies to the left hand side of \eqref{bound} as well. \\
Now, instead of \eqref{narrow}, we assume that the inequality 
$ |v_1 -v_2| \ge 1/N $
holds. Then, by the Christoffel-Darboux formula \eqref{cdsc}, we have 
\begin{equation} \notag
\begin{split}
&|N^{-1}\tilde{K}_N
\bigl(v_1,v_2, N^{-1})|=  N^{1/2}|K_N\bigl(\sqrt{N}v_1,\sqrt{N}v_2)| \\
 & \le \sqrt{N}\bigl( |\varphi_N(\sqrt{N}v_1)| |\varphi_{N-1}(\sqrt{N}v_2)| + 
|\varphi_N(\sqrt{N}v_2)| |\varphi_{N-1}(\sqrt{N}v_1)|\bigr).
\end{split}
\end{equation}
Applying Lemma \ref{hermfunnear} and  enlarging, if necessary, the constant in the bound, the Lemma is proved. 

\end{proof}
\begin{lemma} \label{mixed}
For every $\a \in (0,1)$ and $A>0$ there exists a constant $R_2(\a, A)$ 
such that for every $u \in [-2+\a, 2-\a],$ $t_1, t_2 \in [-A,A],$ $H \in \R$ and 
$N \in \N$ 
we have the relation 
\begin{equation} \label{mix}
|N^{-1}\tilde{K}_N \bigl((u+t_1/N)d(H),(u+t_2/N)d(H), N^{-1})|  
\le \G(N)\bigl(R_2(\a,A)\bigr)^N H^{N}
\end{equation}
whenever none of the following pairs of inequalities 
\begin{equation} \label{pair}
\begin{split}
&\bigl|\Im\bigl({(u+t_i/N)d(H)}\bigr)\bigr| \le \b(\a/4), \quad i=1,2, \\
&\bigl|\Im\bigl({(u+t_i/N)d(H)}\bigr)\bigr| \ge \b(\a/4)/2, \quad i=1,2,
\end{split}
\end{equation}
holds. 
\end{lemma}
\begin{proof}
Set $u_i=(u+t_i/N), v_i=u_id(H), i=1,2.$ Since \eqref{pair} 
does not hold we have 
$|\Im(v_i)| < \b(\a/4)/2$ and   $|\Im(v_j)| > \b(\a/4)$ for some 
choice of $i, j \in \{1,2\},$ $ i \neq j.$ This implies 
\[ \notag |v_2 -v_1| \ge \b(\a/4)/2, \]
 and by the Christoffel-Darboux formula \eqref{cdsc} we have 
\begin{equation} \notag
\begin{split}
&\frac{1}{N}\bigl|\tilde{K}_N
\bigl(v_1,v_2, N^{-1}\bigr)\bigr|= \sqrt{N}\bigl|K_N\bigl(\sqrt{N}v_1,\sqrt{N}v_2\bigr)\bigr| \\
 &\le  \frac{2\sqrt{N}}{\b(\a/4)}\bigl( |\varphi_N(\sqrt{N}v_1)| |\varphi_{N-1}(\sqrt{N}v_2)| + 
|\varphi_N(\sqrt{N}v_2)| |\varphi_{N-1}(\sqrt{N}v_1)|\bigr).
\end{split}
\end{equation}
Assume now that, for example, $|\Im(v_1)| < \b(\a/4)/2$ and   $|\Im(v_2)| > \b(\a/4).$ In addition we have that $ u_1 \in [-2+\a/2,2-\a/2]$ for $N \ge N(\a,A).$ Then, by assertion 2 of Lemma \ref{switchone}, $v_1\in \overline{S}_{\a/4,\b(\a/4)/2} \subseteq \overline{S}_{\a/4,\b(\a/4)},$ and we may apply Lemma  \ref{hermfunnear} (with the parameter value $\a/4$ instead of $\a$) to bound $\varphi_N(\sqrt{N}v_1)$ 
and $\varphi_{N-1}(\sqrt{N}v_1).$ Furthermore, we apply Corollary \ref{pathfar} with $\d = \b(\a/4)$ to bound $\varphi_N(\sqrt{N}v_2)$ and $\varphi_{N-1}(\sqrt{N}v_2).$ 
Therefore, we arrive at the inequality 
\begin{equation} \notag
\begin{split}
\frac{1}{N}|\tilde{K}_N& \bigl((u+t_1/N)d(H),(u+t_2/N)d(H),1/N\bigr)| \\
&\le C(\a,A)\sqrt{N}\G(N) M^{N}(\b(\a/4)) H^{N}.
\end{split}
\end{equation}
Replacing the quantity  $M(\b(\a/4))$ by an appropriate larger constant, we reduce this inequality to the form \eqref{mix}. 
Moreover, this way the bound can be made valid for all $N \in \N,$ which proves the Lemma.

\end{proof}
\begin{lemma} \label{kernelfar} For every $\d>0$ and certain
 constants $C(\d), M(\d)$ the inequality
\begin{equation} \label{kerfar}
\frac{1}{N}|\tilde{K}_N \bigl(z_1,z_2,\tfrac{1}{N})|  \le 
C(\d) \G(N)R^{2N}(\d)\, (\Im(z_1))^N(\Im(z_2))^ N 
\end{equation}
holds for every complex numbers $z_1, z_2$ such that $|\Im z_i| \ge \d$ and $|\Re z_i|
\ge |\Im z_i|$\linebreak 
$(i=1,2),$ and every natural number $N \ge \d^{-2}.$
\end{lemma}
\begin{proof}
First we consider the case when $\Re z_1 \Re z_2 > 0.$ 
Let
$$ f_{a,p}(x)
= \frac{1}{\Gamma(p)} a^{p}x^{p-1} e^{-ax}, p > 0, x> 0,$$ 
be the
{\em gamma density} 
\cite{Fe71}  with the parameters $ a> 0,$ $p> 0,$  
and let $P \ge 0$ and  $N \ge 1$ be some integers. Then for every complex 
number $z$ satisfying $\Re z\ge |\Im z |>0$ we have
\begin{equation} \label{single}
\begin{split}
&\int_0^{\infty}|z+\th|^{2P} e^{-N\Re(z+\th)^2/2}d \th 
=\int_{\Re z}^{\infty}|\xi + i \Im z |^{2P} e^{-\frac{N}{2}\Re(\xi + i
\Im z)^2}d \xi \\
&\le \int_{|\Im z|}^{\infty}|\xi + i \Im z |^{2P}
e^{-\frac{N}{2}\Re(\xi + i \Im z)^2}d \xi 
=\int_{|\Im z|}^{\infty}(\xi^2 +  (\Im z)^2 )^P e^{-\frac{N}{2}(\xi^2 -
(\Im z)^2)}d \xi \\
&= e^{N(\Im z)^2}\int_{|\Im z|}^{\infty}(\xi^2 +
(\Im z)^2 )^P
e^{-\frac{N}{2}(\xi^2 + (\Im z)^2)}d \xi \\
& \le \frac{2^P\, e^{N(\Im z)^2}}{N^{P+1}|\Im z|}\int_{|\Im z|}^{\infty}\bigl(N \bigl(\xi^2
 +  (\Im z)^2\bigr)/2\bigr )^P e^{-\frac{N}{2}(\xi^2
+ (\Im z)^2)}d\bigl(N \bigl(\xi^2+(\Im z)^2\bigr)/2\bigr)\\
&=\frac{2^P\, e^{N(\Im z)^2}}{N^{P+1}|\Im z|} \int_{N(\Im
z)^2}^{\infty} \eta^P
e^{-\eta} d\eta 
=\frac{2^P\, e^{N(\Im z)^2}}{N^{P+1}|\Im z|}\G(P+1)
\int_{N(\Im z)^2}^{\infty}
f_{1,P}(\eta) d \eta \\
&=\frac{2^P\, e^{N(\Im z)^2}}{N^{P+1}|\Im z|}\G(P+1)e^{-N(\Im z)^2}\biggl(1+\frac{N(\Im z)^2}{1!}+\cdots+\frac{[N(\Im z)^2]^P}{P!}
\biggr)\\
&=\frac{2^P\G(P+1)}{N^{P+1}|\Im z|}
\biggl(1+\frac{N(\Im z)^2}{1!}+\cdots +\frac{[N(\Im z)^2]^P}{P!}\biggr),\end{split}
\end{equation}
where we used a well known formula (\cite{Fe71}, p.\! 11)  while
integrating $f_{1,P}.$ 
If $|\Im z| \ge \d$ then $N(\Im z)^2 \ge 1$ for $N \ge \d^{-2},$
and the bound \eqref{single} yields
\begin{equation} \label{farsingle}
\int_0^{\infty}|z+\th|^{2P} e^{-N\Re(z+\th)^2/2}d \th \ \ \le e
2^PN^{-1}\G(P+1)|\Im(z)|^{2P-1}
\end{equation}
whenever $\Re z\ge |\Im z |>0.$ 
In particular, for $P= N$ and $P=N-1$ we have
\begin{equation} \label{farN}
\int_0^{\infty}|z+\th|^{2N} e^{-N\Re(z+\th)^2/2}d \th \ \ 
\le e 2^N N^{-1} \G(N+1)|\Im(z)|^{2N-1}
\end{equation}
and
\begin{equation} \label{farN-1}
\int_0^{\infty}|z+\th|^{2(N-1)} e^{-N\Re(z+\th)^2/2}d \th \ \
\le e 2^{N-1}N^{-1}\G(N)|\Im(z)|^{2N-3}.
\end{equation}
Thus, in view of representation \eqref{Forr} and the bounds in Lemma \ref{hermfunfar}, 
for $z_1, z_2$ such that $ |\Im z_i| \ge \d$ and $\Re z_i\ge|\Im z_i|$ $(i=1,2)$ we have
\begin{equation}
\begin{split}
&|N^{-1}\tilde{K}_N
\bigl(z_1,z_2, N^{-1})| =|N^{-1/2}K_N\bigl(\sqrt{N}z_1,
\sqrt{N}z_2)| \\
&=(2N)^{-1/2}\biggl|\int_0^{\infty}\bigl(
\varphi_N(\sqrt{N}z_1 +\ta)\varphi_{N-1}(\sqrt{N}z_2+\ta) \\
&\hskip3.5cm+\varphi_{N-1}(\sqrt{N}z_1+\ta)
\varphi_N(\sqrt{N}z_2+\ta)\bigr)d\ta \biggr| \\
&=2^{-1/2}\biggl|\int_0^{\infty}\bigl( \varphi_N(\sqrt{N}(z_1
+\th))
\varphi_{N-1}(\sqrt{N}(z_2+\th)) \\
&\hskip3.5cm+\varphi_{N-1}(\sqrt{N}(z_1+\th))
\varphi_N(\sqrt{N}(z_2+\th))\bigr)d\th\biggr|\\
\le &\, C^2(\d)N^{-1/2}M^{2N-1}(\d)\int_0^{\infty} \bigl(|z_1+\th|^N
e^{-N\Re(z_1+\th)^2/4}
|z_2+\th|^{N-1}e^{-N\Re(z_2+\th)^2/4}\\
&\hskip3.5cm+|z_1+\th|^{N-1}e^{-N\Re(z_1+\th)^2/4}
|z_2+\th|^N e^{-N\Re(z_2+\th)^2/4}\bigr)d\th \\
&\le C^2(\d)N^{-1/2}M^{2N-1}(\d)\times\\
&\biggl[
\biggl(\int_0^{\infty}|z_1+\th|^{2N} e^{-N\Re(z_1+\th)^2/2}d \th
\int_0^{\infty}|z_2+\th|^{2(N-1)}
e^{-N\Re(z_2+\th)^2/2}d\th\biggr)^{1/2} \\
\hskip0.3cm+&\biggl(\int_0^{\infty}|z_1+\th|^{2(N-1)}
e^{-N\Re(z_1+\th)^2/2}d \th 
\int_0^{\infty}\bigl(|z_2+\th|^{2N} e^{-N\Re(z_2+\th)^2/2}d \th\biggr)^{1/2} \biggl]. \end{split}
\end{equation}
Combining this bound with relations \eqref{farN} and \eqref{farN-1}
we find that
 \begin{equation} \label{oneside}
\frac{|\tilde{K}_N \bigl(z_1,z_2, N^{-1})|}{N} \le \frac{e
2^{1/2}C^2(\d)M^{2N-1}(\d)2^N\G(N)|\Im(z_1)|^N|\Im(z_2)|^
N}{\d^{2}\,N^{3/2}}.
\end{equation}
The bound \eqref{oneside} applies as well when $z_1, z_2$ satisfy $ |\Im z_i| \ge \d$ 
and $\Re z_i\le -|\Im z|$ $(i=1,2).$ Indeed, this follows from \eqref{symm}. \\
Finally, we consider $z_1, z_2$ such that 
$|\Re z_i| \ge |\Im z|\ge \d$ $(i=1,2)$ and $\Re z_1 \Re z_2 <0.$
In this case we have 
$$ \frac{1}{|z_2 - z_1|}  \le \frac{1}{|\Re(z_1)| + |\Re(z_2)|} 
\le \frac{1}{|\Im(z_1)| + |\Im(z_2)|}\le (2 \d)^{-1}$$
and, by \eqref{cdsc}, 
\begin{equation} \notag
\begin{split}
&\frac{1}{N}|\tilde{K}_N
\bigl(z_1,z_2, N^{-1})|= \frac{1}{\sqrt{N}}|K_N\bigl(\sqrt{N}z_1,\sqrt{N}z_2)| \\
 &\le  \frac{2}{\d \sqrt{N}}\bigl( |\varphi_N(\sqrt{N}z_1)| |\varphi_{N-1}(\sqrt{N}z_2)| + 
|\varphi_N(\sqrt{N}z_2)| |\varphi_{N-1}(\sqrt{N}z_1)|\bigr).
\end{split}
\end{equation}
Applying Lemma \ref{hermfunfar} we obtain
\begin{equation} \label{otherside}
|N^{-1}\tilde{K}_N
\bigl(z_1,z_2, N^{-1})| \le C(\d)N^{-1}M^{2N-1}(\d)2^N\G(N)\, (\Im(z_1))^N(\Im(z_2))^ N. 
\end{equation}
Setting $R(\d)=M(\d)\sqrt{2}$ and comparing \ref{oneside} and \ref{otherside} with the factor $N^{-1}$ omitted, the proof 
is completed. 

\end{proof}
\begin{corollary} \label{far}
Let $\a \in (0,1) $ and $A $ be positive numbers. Then there exist positive constants 
\, $C(\a, A)$ and  $R_2(\a)$ such that for every $N \in \N$
and real numbers \, $u \in [-2,2],$ $t_1,t_2 \in [-A,A],$ and $H$ we have
\begin{equation} \label{scfar}
\left|\frac{1}{N}\tilde{K}_N\! \left(\!\bigl(u+\tfrac{t_1}{N}
\bigr)d(H),\bigl(u+\tfrac{t_2}{N}\bigr)d(H),\tfrac{1}{N}\!\right)\!\right|  \le \G(N) \bigl(R_3(\a, A)\bigr)^{2N} H^{2N}
\end{equation}
whenever the inequalities $\bigl|\Im\bigl({(u+t_i/N)d(H)}\bigr)\bigr|\ge \ \b(\a/4)/2,$ $ i=1,2,$
hold.
\end{corollary}
\begin{proof}
Set $ v_i=(u+t_i/N)\,d(H)$ for $i=1,2.$ 
From \eqref{cone} we conclude that $|\Re v_i| \ge |\Im v_i|,$ $i=1,2.$ 
It follows from \eqref{expbound} that
\[ \notag |\Im(v_i)|^N  \le  C(A) |H|^N, i=1,2.\]
Then by Lemma \ref{kernelfar} with $\d = \b(\a/4)/2$ we have
\begin{equation} \notag
\begin{split}
\left|\frac{1}{N}\tilde{K}_N\! \left(\!\bigl(u+\tfrac{t_1}{N}
\bigr)d(H),\bigl(u+\tfrac{t_2}{N}\bigr)d(H),\tfrac{1}{N}\!\right)\!\right|\\
  \le  C(\b(\a/4)/2,A) &\G(N)R^{2N}(\b(\a/4)/2) H^{2N}
\end{split}
\end{equation}
for every $ N \ge \bigl(\b(\a/4)/2\bigr)^{-2}.$ Substituting $R$ in the latter inequality by a proper larger constant depending also on $A,$ we obtain \eqref{scfar} which is valid for all $N.$ 

\end{proof}

\begin{proposition} \label{last1}
For every $\a \in (0,1)$ and $A>0$ there exists a constant $R(\a, A)$ 
such that for every $u \in [-2+\a, 2-\a],$ $t_1,\dots, t_n \in [-A,A]$ and $H \in \R$ 
we have the relations 
\begin{equation} \label{fin}
\left|\frac{1}{N}\tilde{K}_N\! \left(\!\bigl(u+\tfrac{t_1}{N}
\bigr)d(H),\bigl(u+\tfrac{t_2}{N}\bigr)d(H),\tfrac{1}{N}\!\right)\!\right|  \le  \G(N)\bigl( R(\a, A)(1+H^2)\bigl)^N 
\end{equation}
and
\begin{equation} \label{detfin}
\begin{split}
\frac{|R_{n,N}^{\textsc{GUE},\tfrac{1}{(1+iH)N}}
\!\!(u+\tfrac{t_1}{N},\dots,u+\tfrac{t_n}{N})|}{N^n}& \\
\le n!\bigl(\G(N)\bigr)^n \bigl(R(\a, A)\bigr)^{nN} &(1+H^2)^{n\bigl(N+\tfrac{1}{4}\bigr)}.
\end{split}
\end{equation}
\end{proposition}
\begin{proof}
For every $\a \in (0,1)$ and $A>0,$ every combination of $t_1\in[-A,A],$ $t_2 \in[-A,A],$  $H \in \R,$ and  $N \in \N$ satisfies the assumptions of at least one of the assertions in Lemma \ref{kernelnear}, Lemma \ref{mixed} and Corollary \ref{far}. Hence, at least one of these upper bounds is valid. Writing $ R(\a, A) =\max \bigl(\bigl(R_1(\a,A)\bigr)^2,R_2(\a,A),\bigl(R_3(\a,A)\bigr)^2\bigr),$
 we obtain  
\begin{equation} \notag
\begin{split}
&\Bigl|\frac{1}{N}\tilde{K}_N  \bigl(\bigl(u+\tfrac{t_1}{N}\bigr)d(H),\bigl(u+\tfrac{t_2}{N}\bigr)d(H),\tfrac{1}{ N}\bigr)\Bigr|  \\
&\le  \max\bigl(\bigl(R_1(\a,A)\bigr)^{2N},\G(N)\bigl(R_2(\a, A)\bigr)^{N} 
H^{N},\G(N)\bigl(R_3(\a, A)\bigr)^{2N} H^{2N} \bigr)  \\
&\le \, \bigl( R(\a, A)\bigr)^N \max \bigl(1, \G(N)H^N, \G(N) H^{2N} \bigr) \\
&< \, \bigl( R(\a, A)\bigr)^N \max \bigl(1, \G(N)(1+H^2)^{N/2}, \G(N)  (1+H^2)^N\bigr)\\
&=\G(N)\bigl( R(\a, A)\bigr)^N(1+H^2)^{N},\\
\end{split}
\end{equation}
and the bound \eqref{fin} is proved. \\
Now we will derive    a bound for the correlation function
$R_{n,N}^{\text{GUE},1/((1+iH)N)}.$
In view of formulas \eqref{gdet} and \eqref{sckernel} we have
\begin{equation} \label{sccorr1}
R_{n,N}^{\text{GUE},{\s}s}(x_1,\dots,x_n)=
\s^{-n/2}R_{n,N}^{\text{GUE},s}(\s^{-1/2}x_1,\dots,\s^{-1/2}x_n)
\end{equation}
It follows from \eqref{dissum} that
\begin{equation} \label{sccorr2}
R_{n,N}^{\text{GUE},{\s}s}(x_1,\dots,x_n)
 = \int_0^{\infty} R_{n,N}^{\text{HSE}, u/N^2} \g_{N^2,{\s}s}(u) du,
\end{equation}
where
\begin{equation} \label{gammasc}
\g_{m,s}(u)=\begin{cases}(2^{m/2}s^{-m/2} \G(m/2))^{-1} u^{(m/2)-1}\exp{(-u/(2s))}, \text{if} \,  u \ge 0, \\
0, \text{if} \, u < 0.
\end{cases}
\end{equation}
Since $(x_1,\dots,x_n) \mapsto R_{n,N}^{\text{GUE},s}(x_1,\dots,x_n)$ is 
an entire function, formula \eqref{sccorr1}
can be used to analytically continue the function \linebreak
$R_{n,N}^{\text{GUE},{\s}s}(x_1,\dots,x_n)$ in the parameter $\s$
to the domain $\Re \s >0 .$ On the other hand, an analytic continuation can be obtained from formulas \eqref{sccorr2} and  \eqref{gammasc} (we omit the check of this fact requiring routine bounds for $\g_{m,s}$ in the parameter $s \in \C$). 
Since these continuations obviously agree, both  \eqref{sccorr1} and \eqref{sccorr2} are valid for complex $s$ with $\Re s>0.$   
In particular, we have the identity 
\begin{equation} \label{corresc}
R_{n,N}^{\text{GUE},1/((1+iH)N)}(x_1,\dots,x_n)
=d^n(H)R_{n,N}^{\text{GUE},1/N}(x_1d(H),\dots,x_n d(H)).
\end{equation}
Along with the determinantal representation \eqref{gdet} (which also 
admits an analytic continuation) and inequality \eqref{fin}, this leads to 
\begin{equation} \label{correst}
\begin{split}
&N^{-n}|R_{n,N}^{\text{GUE},1/((1+iH)N)}(u+t_1/N,\dots,u+t_n/N)|\\
&=N^{-n}|d(H)|^n
R_{n,N}^{\text{GUE},1/N)}((u+t_1/N)d(H),\dots,(u+t_1/N) d(H))|
\\
&\le n!|d(H)|^n\max_{0 \le i < j \le
n}|N^{-1}\tilde{K}((u+t_i/N)d(H),(u+t_j/N)d(H))|^n \\
&\le n!|d(H)|^n \bigl( \G(N)\bigl( R(\a, A)\bigr)^N(1+H^2)^N\bigr)^n\\
&=  n!\bigl(\G(N)\bigr)^n \bigl(R(\a, A)\bigr)^{nN} (1+H^2)^{n(N+1/4)}\\
\end{split}
\end{equation}
since $|d(H)|= (1+H^2)^{1/4}$ by \eqref{scalmod}.

\end{proof}

{\em Proof of Theorem.}
As explained in the end of Section 2 we need to show that
\begin{equation} \label{relat}
\begin{split}
 \frac{1}{(Nw(u))^n}\! \!
\int_{-\infty}^{\infty} \! \!  \! \!\phi_{N^2}(h)
R_{n,N}^{\text{GUE},1/((1-ih
\sqrt{2}/N)N)}\! \biggl(\!u+\! & \frac{t_1}{Nw(u)}, \dots,
 u\! + \! \frac{t_n}{Nw(u)} \biggr)dh  \\
 \underset{ N \to \infty}{\longrightarrow}&\!\frac{1}{\sqrt{2\pi}}\det \biggl( \! \frac{\sin \pi(t_i-t_j)}{ \pi(t_i-t_j)}\biggr)_{\mspace{-4.0 mu} i,j=1}^{\mspace{-4.0 mu}n}  \mspace{-16.0 mu}
\end{split}
\end{equation}
uniformly in $u, t_1, \dots, t_n$ subject to the  conditions
formulated in the
the statement of the Theorem. Throughout the proof we will assume
that these conditions are satisfied and omit the arguments of
$R_{n,N}^{\text{GUE},\cdot}$ whenever that possible. The
existence of the integral in the left hand side of 
relation \eqref{relat} will be a consequence of our estimates. \newline
 Note that
$R_{n,N}^{\text{GUE},s}(-u_1,\dots,-u_n) =
R_{n,N}^{\text{GUE},s}(u_1,\dots,u_n)$ since \linebreak
$K_N(-u_1,-u_2)= K_N(u_1,u_2)$ by \eqref{symm}. In view of this property it
suffices to prove \eqref{relat} for $u \ge  0$ only. Observe that
for real $h$
$$\phi_{N^2}(-h)=\overline{\phi_{N^2}(h)}$$ and
$$ R_{n,N}^{\text{GUE},1/((1+ih\sqrt{2}/N)N)}(u_1, \dots, u_n) =\overline{ R_{n,N}^{\text{GUE},1/((1-ih\sqrt{2}/N)N)}}(u_1, \dots, u_n)$$ for 
$u_1, \dots, u_n \in \R$
which shows that \eqref{relat} is established if we prove
\begin{equation} \label{rel}
\begin{split}
& \frac{1}{(Nw(u))^n}
\int_{0}^{\infty}\Re\bigl(\phi_{N^2}(h)
R_{n,N}^{\text{GUE},1/((1+ih \sqrt{2}/N)N)}\bigr)dh \\
 &\hskip3cm\underset{ N \to \infty}{\longrightarrow} \frac{1}{2\sqrt{2\pi}}\det  \biggl( \frac{\sin \pi(t_i-t_j)}{ \pi(t_i-t_j)}\biggr)_{\mspace{-4.0 mu} i,j=1}^{\mspace{-4.0 mu}n}.
\end{split}
\end{equation}
 It follows from the central limit theorem for densities that
\begin{equation} \label{lim}
\int_{0}^{\infty}\Re \bigl(\phi_{N^2}(h)\bigr)dh
 \underset{ N \to \infty}{\longrightarrow}  \frac{1}{2\sqrt{2\pi}}.
\end{equation}
Here $\phi_{N^2}(\cdot)$ is a prelimiting characteristic
function and $1/2\sqrt{2\pi}$ is half the value of the
limiting standard Gaussian density at $0.$ Thus, to prove
\eqref{rel} it suffices to establish the following relation
\begin{equation} \label{last}
\begin{split}
&\frac{1}{(Nw(u))^n} \int_{0}^{\infty}\Re \bigl( \phi_{N^2}(h)
R_{n,N}^{\text{GUE},1/((1+ih \sqrt{2}/N)N)}\bigr)dh \\
&\hskip3cm-\, \det  \biggl( \frac{\sin \pi(t_i-t_j)}{
\pi(t_i-t_j)}\biggr)_{\mspace{-4.0 mu} i,j=1}^{\mspace{-4.0
mu}n}  \mspace{-2.0 mu}\int_{0}^{\infty}
\Re \bigl(\phi_{N^2}(h)\bigr)dh  \\
&=: I(u,N) - J(N) \underset{N \to \infty}{\longrightarrow} 0
\end{split}
\end{equation}
under the same uniformity constraints as above.
Omitting the arguments $t_1,\dots, t_n,$ we set
$$S_n(u,H)=\det  \biggl (
\frac{\sin \bigl(\pi(t_i-t_j)d(H)w(ud(H))/w(u)\bigr)}
{\pi(t_i-t_j)d(H)}
\biggr)_{\mspace{-4.0 mu}
i,j=1}^{\mspace{-4.0 mu}n}
\mspace{-2.0 mu}$$
and 
$$
S_n =S_n(u,0)=\det \biggl( \frac{\sin
\pi(t_i-t_j)}{ \pi(t_i-t_j)}\biggr)_{\mspace{-4.0 mu}
i,j=1}^{\mspace{-4.0 mu}n}.
$$
For a given  $\a \in (0,1),$  
let $\overline{H}(\a)>0$ be a number
whose existence is guaranteed by Lemma \ref{conv}. Let $\e>0$ be an arbitrary
positive number, and $H_0=H_0(\a,A,\e) \in (0,\overline{H}(\a))$ be small enough
to ensure that
\begin{equation} \label{uniconv}
\left|S_n(u,H)-S_n\right| \le
\e2\sqrt{2}/  \pi
\end{equation}
for all $u \in [-2+\a, 2 - \a]$, $|H|\le H_0$ and  $ t_1,
\dots, t_n$ satisfying $|t_1| \le A, \dots,|t_n| \le A.$
Such a number exists because  of continuity of the
$\sin$-kernel and continuity near $0$ of the function
$d(\cdot).$ 
Performing  the change   of
variable $H=h \sqrt{2}/N,$ we may write 
\begin{equation} \label{splitI}
\begin{split}
I(u,N)= &\frac{1}{(Nw(u))^n} \int_{0}^{\infty} \Re \bigl(\phi_{N^2}(h) 
R_{n,N}^{\text{GUE},1/((1+ih \sqrt{2}/N)N)}\bigr) dh  \\
&=  \frac{N}{\sqrt{2}(Nw(u))^n} \int_{0}^{\infty}\Re
\bigl(\phi_{N^2}(NH/\sqrt{2})
R_{n,N}^{\text{GUE},1/((1+iH)N)}\bigr)dH \\
&=  \frac{N}{\sqrt{2}(Nw(u))^n} \int_{0}^{H_0}\dots\, dH +  \frac{N}{\sqrt{2}(Nw(u))^n}\textbf{}\int_{H_0}^{\infty}\dots\, dH \\
& =: I_1(\e,u,N) +I_2(\e,u,N).
\end{split}
\end{equation}
Quite analogously, we have
\begin{equation} \label{splitJ}
\begin{split}
J(N)=&S_n\,\int_{0}^{\infty} \Re \bigl(\phi_{N^2}(h)\bigr)dh = \frac{N \,S_n}{\sqrt{2}} \int_{0}^{\infty}\Re \bigl(\phi_{N^2}(NH/\sqrt{2})\bigr) dH \\
= &\frac{N\,S_n}{\sqrt{2}} \int_{0}^{H_0} \Re \bigl(\phi_{N^2}(NH/\sqrt{2})\bigr) dH + \frac{N\,S_n}{\sqrt{2}}\int_{H_0}^{\infty}\Re \bigl(\phi_{N^2}(NH/\sqrt{2})\bigr) dH \\
=&: J_1(\e,N) + J_2(\e,N).
\end{split}
\end{equation}
 With this notation we need to show that
\begin{equation} \label{diffzero}
I_1(\e,u,N) - J_1(\e,N) \underset{N \to \infty}{\longrightarrow} 0,
\end{equation}
\begin{equation} \label{tail1}
I_2(\e,u,N) \underset{N \to \infty}{\longrightarrow} 0,
\end{equation}
and
\begin{equation} \label{tail2}
J_2(\e,N) \underset{N \to \infty}{\longrightarrow} 0
\end{equation}
uniformly in $u \in [-2+\a, 2-\a]$ and $ t_i \in [-A,A], i=1,\dots,n.$ \\

By Proposition \ref{last1}, for every $\e ,$ we obtain 
\begin{equation} \label{I2}
\begin{split}
&|I_2(\e,u,N)| \\
  &\le\frac{N}{\sqrt{2}(w(u))^n}  n!\bigl(\G(N)\bigr)^n \bigl(R(\a, A)\bigr)^{nN} \int_{H_0}^{\infty }(1+H^2)^{n(N+1/4)}\bigl|\phi_{N^2}(\tfrac{NH}{\sqrt{2}})\bigr| dH. \\
\end{split}
\end{equation}
Since
\begin{equation} \notag
|\phi_{N^2}(NH/\sqrt{2})|= (1+H^2)^{-(N^2/4)}
\end{equation}
and 
\begin{equation} \label{intbound}
\int_L^{\infty}(1+H^2)^{- K}dH 
<\frac{1}{(1+L^2)^{K-1}}\int_0^{\infty}(1+H^2)dH 
= \frac{\pi}{2(1+L^2)^{K -1}}
\end{equation}
for every $K >1,$
we have, in particular, that
\begin{equation} \label{phiint}
\int_L^{\infty}|\phi_{N^2}(NH/\sqrt{2})|dH <
\frac{\pi}{2(1+L^2)^{N^2/4 -1}}.
\end{equation}
Note that $N^2/4 -nN -1/4 >1$ whenever $N >N_0 =2n +\sqrt{4n^2 + 5/4}.$
For every $N > N_0$
we obtain from \eqref{intbound}  with $K=N^2/4 -nN -1/4$ and $L=H_0$ that
\begin{equation} \label{II2}
I_2(\e,u,N) \le  \frac{\pi N}{2\sqrt{2}(w(u))^n} n!\bigl(\G(N)\bigr)^n \bigl(R(\a, A)\bigr)^{nN} 
\frac{ 1}{\bigl(1+H^2_0\bigr)^{N^2/4 -nN -5/4}}.
\end{equation}
Since
\begin{equation} \label{lowdens}
 \inf_{u \in [-2+\e, 2-\e]} |w(u)| = w(2-\e)
\end{equation}
and, by the Stirling formula, $(1+H_0^2)^{N^2/4}$ tends to infinity more rapidly than any fixed power of 
$\G(N)$ (or, moreover, the $N-$th power of any positive number),  we conclude  from
\eqref{II2} that 
\begin{equation} \label{fI2}
I_2(\e,u,N) \underset{N \to \infty}{\longrightarrow} 0
\end{equation}
uniformly in $u \in [-2+\a, 2-\a]$ and $ t_i \in [-A,A], i=1,\dots,n.$ \\
The check of \eqref{tail2} is even more straightforward.  Indeed, using \eqref{intbound} with $K=N^2/4$ and $L= H_0,$
we have

\[ \notag 
\begin{split}
|J_2&(\e,N)| \\
=  &\frac{N\, S_n}{\sqrt{2}}\biggl|\int_{H_0}^{\infty} \Re \bigl(\phi_{N^2}(NH/\sqrt{2})\bigr) dH\biggr| \le\frac{N\, S_n}{\sqrt{2}} \int_{H_0}^{\infty} \bigl|\phi_{N^2}(NH/\sqrt{2})\bigr| dH \\
&= \frac{N\, S_n}{\sqrt{2}} \int_{H_0}^{\infty}  (1+H^2)^{-(N^2/4)} dH =\frac{\pi N\, S_n}{2\sqrt{2}(1+L^2)^{N^2/4 -1}}, 
\end{split}
\]
and \eqref{tail2} follows. \\
Now we  complete the proof by establishing \eqref{diffzero}. 
Recall that we omit the arguments of $R_{n,N}^{\text{GUE},1/N}.$ All
bounds and limit transitions are uniform with respect to $u
\in [-2+ \a,2-\a] $ and $t_1,t_2, \dots, t_n$ satisfying $|t_1|
\le  A, \dots, |t_n| \le  A.$  We have
\begin{equation} \label{null}
\begin{split}
|&I_1(\e,u,N) - J_1(\e,N)|\\ 
&\le \frac{N}{\sqrt{2}}\int_{0}^{H_0}\left|
\left(\frac{R_{n,N}^{\text{GUE},1/((1+iH)N)}} {(Nw(u))^n}- S_n\right)\phi_{N^2}(NH/\sqrt{2})\right| dH \\ 
&\le \frac{N}{\sqrt{2}}\int_{0}^{H_0}\left|\left( \frac{R_{n,N}^{\text{GUE},1/((1+iH)N)}}{(Nw(u))^n}
- S_n(u,H)\right)\phi_{N^2}(NH/\sqrt{2})\right| dH \\
&\hskip1cm+ \frac{N}{\sqrt{2}}\int_{0}^{H_0}\left|\left(S_n(u,H)- S_n\right)\phi_{N^2}(NH/\sqrt{2})\right| dH \\
&=: D^{(1)}_1(\e,u,N)+ D^{(2)}_2(\e,u,N).
\end{split}
\end{equation}
It follows from the determinantal formula \eqref{gdet} and Lemma \ref{conv} that
\begin{equation}
\label{conbou}
 \sup_{0\le H \le H_0}
 \left|\frac{R_{n,N}^{\text{GUE},1/((1+iH)N)}}{(Nw(u))^n}- S_n(u,H)\right|
 \underset{N \to \infty}{\longrightarrow} 0
\end{equation}
uniformly in 
$ u \in [-2+\a,2-\a],$ and $|t_i| \le A \, (i =1, \dots,n).$ 

Therefore, using \eqref{phiint}, we obtain
\begin{equation}
\notag
\begin{split}
 D^{(1)}_1(\e,&u,N)\\
\le \sup_{0\le H \le H_0}
\biggl| \frac{R_{n,N}^{\text{GUE},1/((1+iH)N)}}{(Nw(u))^n}
& - S_n(u,H)\biggr|\, \int_0^{\infty}\bigl|\phi_{N^2}(h)\bigr|dh \underset{N \to \infty}{\longrightarrow} 0
\end{split}
\end{equation} 
since 
\begin{equation} \label{tend} 
\int_0^{\infty}\bigl|\phi_{N^2}(h)\bigr|dh \underset{N \to \infty}{\longrightarrow} \frac{1}{2\sqrt{2 \pi}}.
\end{equation}
Finally, we see from \eqref{uniconv} and  \eqref{tend} that
\begin{equation}
\notag
D^{(2)}_1(\e,u,N) \le \e.
\end{equation}
Since $\e$ is arbitrary small, this completes the proof.
\qed
\vskip 0.7cm
%
%

\end{document}